\documentclass[preprint,11pt]{elsarticle}
\usepackage[utf8]{inputenc}
\usepackage[scale=0.8]{geometry}
\usepackage{amssymb}
\usepackage{amsthm}
\usepackage{mathtools}
\usepackage{bm}
\usepackage{xcolor}
\usepackage{algorithm}
\usepackage{algpseudocode}
\usepackage{multirow}
\usepackage{graphicx}
\usepackage{booktabs}
\usepackage{makecell}
\usepackage{cases}
\usepackage{enumitem}
\usepackage{subcaption}
\usepackage{diagbox}
\usepackage{siunitx}

\usepackage{pgfplots}
\pgfplotsset{compat=1.18}
\usetikzlibrary{calc}
\usepgfplotslibrary{statistics}

\definecolor{graph_1}{RGB}{117,112,179}
\definecolor{graph_2}{RGB}{217,95,2}
\definecolor{graph_3}{RGB}{27,158,119}
\definecolor{graph_4}{RGB}{231,41,138}
\colorlet{allreviews}{graph_1!25!graph_2!50!graph_3!75!graph_4}

\newcommand{\ra}{\color{black}}
\newcommand{\rb}{\color{black}}
\newcommand{\rc}{\color{black}}
\newcommand{\rd}{\color{black}}
\newcommand{\rall}{\color{black}}

\newlength{\plotwidth}
\newlength{\plotheight}

\usepackage{caption}
\usepackage[colorlinks=true, allcolors=blue]{hyperref}
\hypersetup{pdfauthor={J. Aghili, E. Franck, R. Hild, V. Michel-Dansac, V. Vigon}}

\numberwithin{equation}{section}

\algrenewcommand\algorithmicrequire{\textbf{Input:}}
\algrenewcommand\algorithmicensure{\textbf{Output:}}

\newcommand{\R}{\mathbb{R}}
\newcommand{\bs}{\boldsymbol}
\DeclareMathOperator*{\argmin}{argmin}

\usepackage[capitalize, nameinlink, noabbrev]{cleveref}

\makeatletter
\def\ps@pprintTitle{%
    \let\@oddhead\@empty
    \let\@evenhead\@empty
    \def\@oddfoot{\footnotesize\itshape
         {Submitted preprint} \hfill\today}%
    \let\@evenfoot\@oddfoot
    }
\makeatother

\begin{document}

\begin{frontmatter}

    \title{Accelerating the convergence of Newton's method for nonlinear elliptic PDEs using Fourier neural operators}

    \author[Inria,IRMA]{Joubine Aghili}
    \ead{aghili@unistra.fr}

    \author[Inria]{Emmanuel Franck}
    \ead{emmanuel.franck@inria.fr}

    \author[GE]{Romain Hild}
    \ead{hild.romain@gmail.com}

    \author[Inria]{Victor Michel-Dansac\corref{cor1}}
    \ead{victor.michel-dansac@inria.fr}
    \cortext[cor1]{corresponding author}

    \author[Inria,IRMA]{Vincent Vigon}
    \ead{vigon@math.unistra.fr}

    \affiliation[Inria]{organization={Université de Strasbourg, CNRS, Inria, IRMA},
        postcode={F-67000},
        city={Strasbourg},
        country={France}
    }

    \affiliation[IRMA]{organization={IRMA, Université de Strasbourg, CNRS UMR 7501},
        addressline={7 rue René Descartes},
        postcode={67084},
        city={Strasbourg},
        country={France}}

    \affiliation[GE]{organization={GE Healthcare},
        addressline={2 Rue Marie Hamm},
        postcode={67000},
        city={Strasbourg},
        country={France}}

    \begin{abstract}
        It is well known that Newton's method
        can have trouble converging if
        the initial guess is too far from the solution.
        Such a problem particularly occurs when this method is used
        to solve nonlinear elliptic partial differential equations (PDEs)
        discretized via finite differences.
        This work focuses on accelerating Newton's method convergence in this context.
        We seek to construct a mapping
        from the parameters of the nonlinear PDE
        to an approximation of its discrete solution,
        independently of the mesh resolution.
        This approximation is then used as an initial guess for Newton's method.
        To achieve these objectives, we elect to use a Fourier neural operator (FNO).
        The loss function is the sum of a data term
        (i.e., the comparison between known solutions and outputs of the FNO)
        and a physical term (i.e., the residual of the PDE discretization).
        Numerical results, in one and two dimensions,
        show that the proposed initial guess accelerates the convergence of Newton's method
        by a large margin compared to a naive initial guess,
        especially for highly nonlinear and anisotropic problems, with larger gains on coarse grids.
    \end{abstract}


\end{frontmatter}

\section{Introduction}

The broad context of this work is the resolution of systems of nonlinear equations
using Newton's method.
Many applications (nonlinear elliptic partial differential equations (PDEs),
implicit time-stepping for nonlinear diffusion or hyperbolic equations, \dots)
require inverting large nonlinear systems.
Throughout this document, such systems will be denoted by
\begin{equation}
    \label{eq:main_system_intro}
    F(u) = 0,
\end{equation}
where $F: \mathbb{R}^{N_h} \to \mathbb{R}^{N_h}$ is a known nonlinear function
and $u \in \mathbb{R}^{N_h}$ is the unknown vector.
The integer $N_h$ represents the size of the vector $u$;
in the case of a nonlinear system arising from the discretization of a PDE,
$N_h$ would be the number of degrees of freedom of the discretization.

For smaller problems, Newton's method is often used.
It consists in linearizing the system
around a known state $u_0 \in \mathbb{R}^{N_h}$, to obtain
an affine approximation of $F$ in a neighborhood of $u_0$:
\[
    F(u_0)+F'(u_0) (u-u_0) \approx 0.
\]
This equation is nothing but a linear system with unknown $u$,
whose matrix is the Jacobian matrix of $F$.
To apply Newton's method, this linear system is then repeatedly solved
using the previous iteration as a reference state.
This leads to the following iterative process,
where $u^{(k)}$ is expected to tend towards~$u$ as $k$ goes to infinity.
\begin{algorithm}
    \begin{enumerate}[nosep]
        \item[(1.)] \label{item:initialization_step}
              initialization step: set $u^{(0)} = u_0$,
        \item[(2.)] main loop: for $k \geq 0$,
              \begin{enumerate}[nosep]
                  \item[(2.a.)] \label{item:solve_step}
                        solve the linear system $F(u^{(k)}) + F'(u^{(k)}) \delta^{(k+1)}=0$
                        for $\delta^{(k+1)}$,
                  \item[(2.b.)] \label{item:update_step}
                        update $u^{(k+1)} = u^{(k)} + \delta^{(k+1)}$.
              \end{enumerate}
    \end{enumerate}
    \caption{Newton's algorithm to solve \eqref{eq:main_system_intro}.}
    \label{algo:Newton}
\end{algorithm}

In practice, the different steps of \cref{algo:Newton}
have to be adapted to ensure convergence.
For instance, the algorithm has to be properly initialized,
which means that the initial guess $u_0$
in step~\href{item:initialization_step}{(1.)}
has to be close enough to the solution.
This is the focus of the present work:
devising a novel, data-driven method
to predict a suitable initial guess $u_0$.
Moreover, for larger systems,
such as the ones stemming from the discretization of a PDE,
it is usual to use Jacobian-Free Newton-Krylov (JFNK)
methods \cite{knoll2004jacobian} or
inexact Newton methods \cite{ern2013adaptive}.
In these methods, the linear systems are solved using an iterative method
like the generalized minimal residual method (GMRES) or
the conjugate gradient method (CG).
This means replacing step \href{item:solve_step}{(2.a.)} of \cref{algo:Newton}
with an iterative linear solver.
Such iterative methods only compute matrix-vector products with the Jacobian matrix.
Consequently, it is not necessary to construct and store the Jacobian matrix,
but it is sufficient to use an approximation of the matrix-vector product,
defined for $v \in R^{N_h}$ by
\[
    F'(u)v\approx\frac{F(u+\varepsilon v)-F(u)}{\varepsilon},
\]
and $\varepsilon$ is correctly chosen.
When the linearization error is large,
the linear solver does not need to converge to machine precision.
Therefore, in the JFNK method,
the linear solver threshold is adapted to the nonlinear residual.
In practice, the linear solver is stopped as soon as
the following criterion is satisfied:
\[
    \| F(u^{(k)})+ F'(u^{(k)}) \delta^{(k+1)}\| \leq \eta^{(k)} \| F(u^{(k)})\|,
\]
where the threshold $\eta^{(k)}$ depends on the nonlinear convergence,
see~\cite{osti_218521,2007JCoAM.200,RePEc:spr:annopr}.

This JFNK method is, however, not always sufficient.
Indeed, if the problem is strongly nonlinear,
the method can converge very slowly, and even sometimes fail to converge.
Additionally, for multiscale PDEs like
the magnetohydrodynamics (MHD) equations \cite{franck2015energy}
or anisotropic elliptic equations \cite{degond2012asymptotic,deluzet2019two},
the Jacobian matrix can be ill-conditioned,
which makes the linear step difficult to solve.
For this second point,
linear preconditioning \cite{CheChaLeiKnoTai2014} can be applied to the method,
but could fail to improve the nonlinear convergence.
There exist several techniques to accelerate this convergence;
we provide a brief state of the art of such approaches
in \cref{sec:accelerating_newton}.
In this work, we develop a data-driven method to
accelerate the convergence of Newton's method by
enhancing the initialization
step~\href{item:initialization_step}{(1.)}
of \cref{algo:Newton}.
This methodology, as well as the problems we apply it to,
are described in \cref{sec:context}.
{\rall Finally, the outline of the paper is given in \cref{sec:outline}.}

\subsection{Convergence acceleration: brief overview of existing approaches}
\label{sec:accelerating_newton}

In this paragraph, we give a brief overview of several techniques
dedicated to the acceleration of Newton's method.
We refer the reader to \cite{knoll2004jacobian} for a more complete review.
Namely, we present line search in \cref{sec:line_search},
adaptive inexact Newton methods in \cref{sec:adaptive_inexact_newton},
nonlinear preconditioning in \cref{sec:nonlinear_preconditioning},
and initial guess prediction in \cref{sec:initial_guess_prediction}.
This paragraph comes in support of the choice we made in
this paper to accelerate Newton's method convergence
through a better prediction of the initial guess.

\subsubsection{Line search}
\label{sec:line_search}

A classical approach to accelerate the convergence of the JFNK method
is to use line search \cite{BroSaa1990,knoll2004jacobian,Pawlowski2006GlobalizationTF}.
This helps to globalize the JFNK method,
i.e. to relax the choice of the initial guess $u_0$.
The idea is to compute a descent direction $\delta^{(k+1)}$
from the classical Newton linear solve,
i.e. step \href{item:solve_step}{(2.a.)} of \cref{algo:Newton}.
Then, at each iteration of the Newton solver,
the linear search consists of finding iteratively a large enough real number $\lambda \leq 1$ such that
\[
    \| F(u^{(k)} + \lambda \delta^{(k+1)})\|
    \leq
    \| F(u^{(k)})\|,
\]
and then replacing the classical Newton update
$u^{(k+1)}=u^{(k)}+\delta^{(k+1)}$
(corresponding to step \href{item:update_step}{(2.b.)} of \cref{algo:Newton})
with $u^{(k+1)}=u^{(k)}+\lambda \delta^{(k+1)}$.
To stop the iterative process of the line search, there exist many criteria,
among which the Armijo or Wolfe conditions \cite{PhilipCONVERGENCECF}
are most often used.

\subsubsection{Adaptive inexact Newton methods}
\label{sec:adaptive_inexact_newton}

In the specific case of nonlinear systems stemming from
a finite element discretization of a PDE,
there exists a strategy using a posteriori error estimates.
This method was first proposed in~\cite{ern2013adaptive}
for nonlinear elliptic equations,
and then extended to multiphase flows in porous media
in~\cite{di2014posteriori,di2014posteriori2}.
First, an a posteriori estimate is constructed
using the finite element method to discretize the PDE.
This estimate is then employed to adaptively change the linear and nonlinear
convergence criteria, while simultaneously refining the mesh.
Using such estimates reduces the global cost of the method for a given accuracy.
While efficient, this class of methods suffers from a lack of generality,
since the a posteriori estimate has to be written
for the problem and the numerical method under consideration.
Additionally, if the initial guess $u_0$ is too far from the exact solution,
then convergence issues may remain.

\subsubsection{Nonlinear preconditioning}
\label{sec:nonlinear_preconditioning}

Another approach is the use of nonlinear preconditioners.
It is an extension of classical linear preconditioning.
Consider a strongly nonlinear system $F(u)=0$.
Nonlinear preconditioning consists of finding a nonlinear operator~$G$
such that $G \circ F \approx A$, with $A$ a linear operator.
In addition, if $G \approx F^{-1}$,
then $A\approx \text{Id}$, which accelerates the convergence.
Applying nonlinear preconditioning boils down to solving
another (almost linear) nonlinear system
\[
    G(F(u))=G(0),
\]
where a fast convergence is expected,
since the operator $G \circ F$ is close to being linear.
The main difficulty of this approach, in addition to finding $G$,
is ensuring that the Jacobian matrix of $G \circ F$ (or an approximation thereof)
is nonsingular and easily computable,
to ensure that the linear solve,
step~\href{item:solve_step}{(2.a.)} of \cref{algo:Newton},
has a satisfactory convergence.

In \cite{cai2002nonlinearly,dolean2016nonlinear},
the authors propose nonlinear Schwarz-based preconditioning
for Newton's method in the context of domain decomposition.
It is an extension of a method commonly used to precondition linear problems.
In \cite{tang2019fully}, the authors design
a physics-based nonlinear preconditioner for
potentially discontinuous, time-independent PDEs.
Their approach is based on detecting and eliminating
strongly nonlinear or discontinuous regions.

\subsubsection{Initial guess prediction}
\label{sec:initial_guess_prediction}

Like any iterative algorithm, Newton's method requires an initial guess $u_0$
close enough to the actual solution $u$.
An alternative approach to accelerate the convergence of the method is
the creation of an initial guess $u_0$,
ensuring that $u_0$ is close enough to $u$.
This acts on the initialization step
\href{item:initialization_step}{(1.)} of \cref{algo:Newton}.
To that end,
an immediate idea is to construct an operator $G$ which approximates $F^{-1}$,
and to use $G(0)$ as an initial guess of Newton's method.
Compared to nonlinear preconditioning, this approach may be more restrictive on $G$
but eliminates the nonsingular requirement on the Jacobian matrix of $G \circ F$.

As explained above, a generic possibility to construct this initial guess is to involve a nonlinear preconditioner.
A second one, for nonlinear systems stemming from PDE discretizations,
is to use discretization- or physics-based criteria.
For time-dependent problems, one can use the previous time step,
which was often (not always) sufficient to
make Newton's method converge (albeit potentially slowly).
One can also use the solution to a linearization of the equation
at the current time step, see~\cite{10.1016/j.camwa.2022.04.013}.
However, these guesses cannot be used for elliptic equations.
Another strategy is to use the solution of a simpler equation as an initial guess.
For example, one can solve a (linear) Stokes problem to obtain an initial guess
for nonlinear PDEs such as the Navier-Stokes or MHD equations,
see~\cite{10.1016/j.cam.2016.06.022}.

Several approaches based on machine learning have also been
proposed to predict an initial guess.
For instance,
see \cite{HuaWanYan2020} for more general iterative methods,
\cite{odot2022deepphysics} for nonlinear elasticity problems,
\cite{LuoCai2023} for applications in high Reynolds number compressible flows
or \cite{NovPoeLugPulCon2024} to accelerate the simulation of chemical reactions
in the context of a coupling with fluid dynamics.

The main idea behind the current work is similar,
but makes use of different neural networks with a more robust learning approach,
as explained in \cref{sec:learning_initial_guess}.

\subsection{Target problem and chosen methodology}
\label{sec:context}

{\rall
    In this work, we propose to accelerate the convergence of the JFNK method
    by constructing a suitable initial guess.
    To that end,
    we present the target problem (nonlinear elliptic PDEs) in \cref{sec:target}
    and the methodology to construct an initial guess in \cref{sec:methodology}.

    \subsubsection{Target problem}
    \label{sec:target}

    This work focuses on elliptic problems
    since there usually is no information on how to choose a suitable initial guess.
    Indeed, it is often harder to obtain a suitable initial guess
    for nonlinear elliptic PDEs compared to time-dependent PDEs.
    The main reason for this is that there is no time stepping,
    and thus the solution at the previous time step
    can no longer be used as an initial guess.
    Note that the proposed method, however, can also be used
    to improve the initial guess for time-dependent problems.
}

For the sake of simplicity in the notation,
all continuous quantities are denoted with capital letters,
and discrete ones with lowercase letters;
all continuous operators are denoted with calligraphic letters,
and discrete ones with capital letters.

We consider both
one-dimensional (1D) and two-dimensional (2D) elliptic problems,
described below.
They are governed by similar equations,
but with slightly different parameters,
and thus we write both problems separately.

The 1D problem is governed by the
following elliptic equation on $\Omega\subset\mathbb{R}$,
with unknown $U: \Omega \to \mathbb{R}$:
\begin{equation}
    \label{1deq}
    \begin{dcases}
        U(x) - \alpha_0 \, \partial_x\left(
        K(x) \, | U(x)|^p \, \partial_x U(x)
        \right) = \Phi(x),
         & \quad \forall x\in \Omega,         \\
        U(x)=0,
         & \quad \forall x\in \partial\Omega,
    \end{dcases}
\end{equation}
with $\Phi \in L^2(\Omega)$
and $K \in \mathcal{C}^0(\Omega, \mathbb{R})$ two given functions,
$p$ even, and $\alpha_0 > 0$.
We assume that there exists $K_0>0$ such that, for all $x\in\mathbb{R}$,
$K(x) > K_0 (1 + x^2)$.
For simplicity,
homogeneous Dirichlet boundary conditions are prescribed on $\partial\Omega$.

The 2D problem is a
general anisotropic nonlinear diffusion equation
on $\Omega\subset\mathbb{R}^2$,
with unknown $U: \Omega \to \mathbb{R}$ and governed by:
\begin{equation}
    \label{2deq}
    \begin{dcases}
        U(x) - \nabla \cdot \left(
        K(x) \, | U(x)|^p \, \nabla U(x)
        \right) = \Phi( x),
         & \quad \forall  x\in \Omega,         \\
        U(x)=0,
         & \quad \forall  x\in \partial\Omega,
    \end{dcases}
\end{equation}
{\rb
with $p$ an even number,
and where $\Phi \in L^2(\Omega)$
and $K \in \mathcal{C}^0(\Omega, \mathcal{M}_2(\mathbb{R}))$
are two given functions.
We have denoted by $\mathcal{M}_2(\mathbb{R})$ the space of
$2 \times 2$ real matrices.}
Homogeneous Dirichlet boundary conditions
are once again prescribed on $\partial\Omega$.
We assume $K$ is bounded by below, i.e. there exists $K_0>0$ such that, for all ${x}\in\mathbb{R}^2$,
$\langle {x},K({x}){x} \rangle > K_0 | {x}|^2$.
Such problems are common in physical applications;
it appears for example in tokamak simulations \cite{hoelzl2021jorek}
or simulations of multiphase flows in porous media \cite{AghBreHenMasTre2019}.
For $p=0$, a classical linear anisotropic diffusion equation is obtained.
In this case, the asymptotic limit is
ill-posed~\cite{degond2012asymptotic,deluzet2019two},
which may highly degrade the conditioning.

Both problems \eqref{1deq} and \eqref{2deq}
can be rewritten in a more compact form.
Indeed, assume that
the solution $U$ belongs to some Hilbert space $\mathcal{U} \subset L^2(\Omega)$,
and define the function
\begin{equation*}
    A
    \, = \,
    (\Phi, K)
    \, \in \,
    \mathcal{A}
    \, \subset \,
    L^2(\Omega) \times \mathcal{C}^0(\Omega, \mathcal{M}_d(\mathbb{R})),
\end{equation*}
where $d$ is the dimension of the space domain $\Omega$
{\rb and $\mathcal{M}_d(\mathbb{R})$ is the space of
        $d \times d$ real matrices.}
Therefore, $A$ corresponds to the data of the physical model.
Therefore, \eqref{1deq} and~\eqref{2deq} can be rewritten as
\begin{equation}
    \label{eq:general_pde}
    \mathcal{F}(U; A) = 0,
\end{equation}
where $\mathcal{F}: \mathcal{U} \times \mathcal{A} \to L^2(\Omega)$
is the residual operator,
with natural boundary conditions for simplicity.
In the remainder of this work, we consider $\mathcal A$
as a set of ``parameters'' of the physical model.
In reality, $\mathcal A$ corresponds to a set of functions
parameterizing the physical model.

    {\rall
        \subsubsection{Main goal}
        \label{sec:goal}

        The broad scope of this work is to
        accelerate high-fidelity PDE simulations,
        i.e. accelerate the resolution of the discrete version of \eqref{eq:general_pde}.
        In this context, a numerical scheme is applied to \eqref{eq:general_pde},
        which provides rigorous convergence guarantees as the mesh is refined.
        For such nonlinear problems, numerical schemes require the use
        of a nonlinear solver like Newton's method.
        For a given mesh,
        we know that Newton's method will converge,
        provided a suitable initial guess is provided.
        Such convergence guarantees are of paramount importance
        in many applications,
        such as nuclear power plant simulations \cite{HelHurQui2020}
        or aeronautics \cite{DecGanBruBen2014}.

        Moreover, note that the source term $\Phi$ and the diffusion matrix $K$
        in \eqref{1deq} and \eqref{2deq} are generic.
        In practice, this is the case when considering the discretization of time-dependent PDEs,
        or when studying the solution with respect to a parameter in uncertainty quantification.
        The goal of this work is therefore to provide a unified procedure for
        a wide range of such source terms and diffusion matrices.
        This will be done by solving a single equation \eqref{eq:general_pde},
        but whose parameters $A$ can vary widely.
        All in all, we seek to obtain a high-fidelity, provably convergent model
        for a wide range of such source terms and diffusion matrices.

        This last remark is why we choose to predict the initial guess.
        Indeed, this will make it possible to benefit from the convergence properties
        of Newton's method.
        Moreover, in general, finding such an initial guess is a difficult task.
        For instance, \cref{sec:motivation} demonstrates the failure of
        several naive (and not so naive) initial guesses.
        We therefore need to construct a good prediction of the initial guess:
        the method we elect to use is described in the next section.
    }

    {\rall
        \subsubsection{Methodology to predict an initial guess}
        \label{sec:methodology}
    }

As mentioned before, we propose to construct a method
able to produce a good initial guess for the JFNK method.
We discretize the generic PDE \eqref{eq:general_pde}
using a classical finite difference scheme with $N_h$ mesh points.
We obtain the following system of nonlinear equations:
\begin{equation}
    \label{eq:discretized_pde}
    F(u; a) = 0,
\end{equation}
with $u \in \R^{N_h}$ the discretization of the unknown function $U$,
and $a$ a discretization of $A$,
i.e. a discretization of $\Phi$ and $K$.
Analogously to $A = \{\Phi, K\}$, we write $a = \{\varphi, k\}$,
with $\varphi \in \R^{N_h}$ a discretization of $\Phi$
and $k \in \R^{n_c-1}$ a discretization of $K$,
where $n_c = 2$ in 1D for $d=1$ and $n_c = 5$ in 2D for $d=2$.
Therefore, $F: \R^{N_h} \times \R^{n_c N_h} \to \R^{N_h}$
is a discretization of the operator $\mathcal{F}$.
Our objective is then to construct an operator
$G^{+}: \R^{n_c N_h} \to \R^{N_h}$
such that
\begin{equation}
    \label{eq:G_theta_approx_F_inv}
    F(G^+(a);a) \approx 0.
\end{equation}

Note that, if \eqref{eq:G_theta_approx_F_inv} were an equality,
$G^+$ would be the pseudo-inverse of
$F$ with respect to its first variable:
we therefore seek to approximate this pseudo-inverse.
This function $G^{+}$ can be viewed as a discrete operator
mapping the data on the mesh to the solution on the mesh,
contrary to the nonlinear preconditioners
mentioned in \cref{sec:nonlinear_preconditioning},
which approximates the inverse of the continuous operator.
Since $G^{+}$ is a function from $\R^{n_c N_h}$ to $\R^{N_h}$,
it quickly becomes high-dimensional when the mesh is refined.

For this reason, we propose to use a neural network
with trainable weights $\theta$ to construct a function~$G_\theta^+$
approximating the inverse of $F$, like in \eqref{eq:G_theta_approx_F_inv}.
Indeed, for the last ten years, neural networks have shown their ability
to outperform other methods for high-dimensional regression problems.
Such operators constructed via neural networks are called neural operators;
the reader is directed to \cite{kovachki2021neural} for
a description of a unified framework for operator learning.

Operator learning comes in two flavors:
the discrete case, where a map between discretizations of functions is built,
and the continuous case, where the map is between the actual functions themselves.
On the one hand, for the discrete case,
the most natural approach is to use a convolutional neural network (CNN)
to construct $G_\theta^+$ on a structured grid.
Such neural networks have been successfully implemented
to learn mappings between functions in various applications,
see e.g.~\cite{bois2020neural,geist2021numerical,maulik2021reduced,10.1016/j.jcp.2021.110928,fuhg2023deep,long2018pde}.
The main drawback of this approach is that the mesh used for the learning process
partially restricts the domain of applicability of the discrete operator.
On the other hand, in the continuous approach,
the operator is actually constructed between functions,
mapping $\Phi$ and $K$ to $U$.
Several approaches have been proposed to construct such operators,
see for instance~\cite{li2020neural,li2020multipole,lu2021learning},
which have been unified in a general framework in \cite{kovachki2021neural}.

The approach best suited to our case seems to be Fourier neural operators (FNOs),
introduced in~\cite{li2020fourier}.
{\ra%
Indeed, when considering classical neural operators such as FNOs or DeepONets,
it seems like FNOs outperform DeepONets in most cases,
see for instance \cite{raonic2023convolutional}, Table 1.%
}
The idea behind FNOs is to learn a convolution in the frequency domain
instead of the usual space domain.
It makes it possible to obtain mesh- and discretization-independent operators.
We will use an FNO to approximate the continuous operator,
before applying any discretization to obtain
an approximation of the pseudo-inverse of $F$.

    {\rd%
        An interesting approach to construct an initial guess
        in the Int-Deep method from~\cite{HuaWanYan2020}.
        Int-Deep uses a partially trained neural network to provide an initial guess.
        As such, one prediction is not very costly,
        but this training has to be restarted for each new problem.
        However, applications such as time-dependent PDEs or uncertainty quantification
        require solving the nonlinear equation many times.
        Each Newton solve may take a substantial amount of time,
        depending on how good the initial guess is.
        In such cases, the training time of the FNO quickly becomes negligible
        if it manages to provide a good initial guess.%
    }

    {\rall%
        Another possibility would be to directly solve the continuous PDE
        \eqref{eq:general_pde} using a neural network~\cite{raissi2019physics, EYu2018}
        or a neural operator~\cite{kovachki2021neural,raonic2023convolutional}.
        However, at the moment, such approaches lack a robust convergence framework,
        contrary to Newton's method.
        Indeed, for a given mesh and a suitable initial guess,
        we know that Newton's method is able to converge.
        Simply using a neural network or a neural operator
        would not provide such theoretical guarantees.
        Another advantage lies in avoiding the strong drift in time
        of the FNO solution compared to the reference
        in complex time-dependent applications such as
        nuclear fusion~\cite{GopPamZanLiGraBreBhaStaKusDeiAna2024}.
    }

\subsection{\texorpdfstring{\rall Outline}{Outline}}
\label{sec:outline}

{\rall%
    At this level, we have defined the continuous problem \eqref{eq:general_pde}
    under consideration, as well as its discrete version \eqref{eq:discretized_pde}.
    Furthermore, we have elected to use an FNO to construct an initial guess.
    The remainder of the paper is dedicated to the learning process,
    and to validate our approach.%
}
To that end, in \cref{sec:learning_initial_guess},
we describe the whole process of learning this initial guess.
Namely, we discuss the choice of the neural network,
the loss functions, and the data generation process.
Then, in \cref{sec:numerical_results}, numerical results
are given in 1D and 2D to validate our approach.
A conclusion ends the paper in \cref{sec:conclusion}.

\section{Learning an initial guess}
\label{sec:learning_initial_guess}

{\rall The goal of this section is to describe the learning process of the initial guess.}
It is organized as follows: firstly, the loss functions are
introduced in \cref{sec:loss_function}. Secondly, we discuss the generation of
training, validation and testing datasets in \cref{sec:data_generation}.
Thirdly, the structure of FNOs is briefly sketched in \cref{sec:NN_structure}.
Fourthly, the hyperparameters of the network are stated, and determined with a grid search, in \cref{sec:hyperparameter_choice}.

\subsection{Discretization-Informed loss function}
\label{sec:loss_function}

Before introducing the full learning process
and the specific network architecture considered in this work,
we start by discussing and motivating the loss functions to be minimized.
Indeed, we have to correctly manage the interplay between
data, physics and discretization.
For the sake of simplicity, we consider the generic PDE \eqref{eq:general_pde},
governed by $\mathcal{F}(U; A) = 0$.

To define our loss functions,
let us start by temporarily remaining at the continuous level.
The best pseudo-inverse operator $\mathcal{G}^+ : \mathcal{A} \to \mathcal{U}$
would satisfy
\begin{equation}
    \label{eq:continuous_l2loss}
    \left \|
    U - \mathcal{G}^+(A)
    \right \|^2_{L^2(\Omega)}
    = 0,
\end{equation}
for all $U$ and $A$ such that $\mathcal{F}(U; A) = 0$.
This is nothing saying that $\mathcal{G}^+(A) = U$ in the $L^2$ sense.
Similarly, we can replace the $L^2$ norm
in the equation above by the $H^1_0$ norm, to obtain
\begin{equation}
    \label{eq:continuous_h1loss}
    \left \|
    \nabla U - \nabla \mathcal{G}^+(A)
    \right \|_{L^2(\Omega)}^2
    = 0.
\end{equation}
In addition, since $\mathcal{G}^+$ is such that
$\mathcal{F}(\mathcal{G}^+(A); A) = 0$, it also satisfies
\begin{equation}
    \label{eq:continuous_physical_loss}
    \left \|
    \mathcal{F}(\mathcal{G}^+(A); A)
    \right \|_{L^2(\Omega)}^2
    = 0.
\end{equation}

Now, we look for an approximation $G_\theta^+$
of this ideal operator $\mathcal{G}^+$
among a given set of parameterized operators.
To that end, we will build it such that it minimizes
a weighted sum of three different discrete loss functions,
consistent with the three equations
\eqref{eq:continuous_l2loss}, \eqref{eq:continuous_h1loss} and
\eqref{eq:continuous_physical_loss}.
Recall that the discretized PDE reads $F(u; a) = 0$.
From now on, we assume that we have at our disposal
$N_\text{data}$ data points $(u_j, a_j)_{j \in \{1, \dots, N_\text{data}\}}$
such that
\begin{equation*}
    \forall j \in \{1, \dots, N_\text{data}\}, \quad
    u_j \in \R^{N_h}, \quad
    a_j \in \R^{n_c N_h}, \quad
    \text{and } F(u_j; a_j) = 0.
\end{equation*}
Generating this data is the object of \cref{sec:data_generation}.

First, we wish to minimize, for all $j$,
the $L^2$ error between the prediction~$G_\theta^+(a_j)$ of the neural network
and the solution $u_j$ of the discretized PDE.
This is the discrete analog of \eqref{eq:continuous_l2loss},
which reads
\begin{equation}
    \label{eq_loss2}
    \mathcal{L}_{\text{data}}^{L^2}(\theta)
    =
    \frac 1 {N_\text{data}}
    \sum_{j=1}^{N_\text{data}} \big\| u_j - G_\theta^+(a_j) \big\|^2.
\end{equation}
This is a classical loss function used in supervised learning.

Second, minimizing the $H^1$ error instead of the $L^2$ error
to provide a discrete analog of \eqref{eq:continuous_h1loss}
requires defining a discrete approximation of the gradient.
For simplicity, we use a centered discretization of the first derivative,
which we denote by $D$.
This leads to the following loss function:
\begin{equation}
    \label{eq_loss3}
    \mathcal{L}_{\text{data}}^{H^1}(\theta)
    =
    \frac 1 {N_\text{data}}
    \sum_{j=1}^{N_\text{data}} \Big\| D \big( u_j - G_\theta^+(a_j) \big) \Big\|^2.
\end{equation}

Finally, what is usually done in the context of
Physics-Informed Neural Operators (PINOs)
is to take the PDE into account in the loss function,
see \cite{WanWanPer2021,li2021physics}.
Therefore, learning such operators
mix classical operator learning and
Physics-Informed Neural Networks
(PINNs)~\cite{raissi2019physics,karniadakis2021physics}.
This would lead to a loss function analog to
\eqref{eq:continuous_physical_loss},
where $G_\theta^+$ would directly approximate $\mathcal{G}^+$.

However, in our case, the prediction $G_\theta^+(A)$ of the neural network
is to be used as an initial guess for Newton's method.
    {\rb%
        Indeed, recall that the goal of this work is to solve nonlinear systems
        stemming from the discretization of (potentially time-dependent) PDEs.
        In this discrete framework,
        we therefore wish to solve the discrete problem, rather than the continuous one.
        In this case, finding $G_\theta^+$ close to~$\mathcal{G}^+$~(in some sense)
        has no guarantee of minimizing the residual of the numerical scheme.%
    }
Indeed, it is entirely possible that $G_\theta^+$,
viewed as an approximation of the continuous operator $\mathcal{G}^+$,
would fail to minimize the discrete residual
$F(G_\theta^+(a_j); a_j)$ for some~$j$,
leading to a poor initial guess for Newton's method.
For this reason, we design another loss function,
based on the PDE discretization $F$,
which we call \emph{discretization-informed}.
This loss function takes the discretization $F$ into account
rather than the continuous PDE $\mathcal{F}$, and is defined by
\begin{equation}
    \label{eq_loss1}
    \mathcal{L}_{\text{dis}}(\theta)
    =
    \frac 1 {N_\text{data}}
    \sum_{j=1}^{N_\text{data}} \Big\| F(G_\theta^+(a_j); a_j) \Big\|^2.
\end{equation}
{\rb%
Note that using the discrete $F$ rather than the continuous $\mathcal{F}$
in this loss function has another advantage:
it is not necessary to compute second-order space derivatives
of the network's output,
which would be cumbersome and computationally expensive.%
}

\subsection{Data generation}
\label{sec:data_generation}

Now that the loss functions have been defined,
we describe how to generate training and validation data.
In this paragraph, we assume that the nonlinearity $p \in \mathbb{N}$ is fixed.
Moreover, we denote by $E$ the discrete elliptic operator,
defined such that $E(u; k) = \varphi$.
Note that $E$ is simply another formulation of~$F$:
indeed, $E(u; k) = \varphi$ is equivalent to $F(u; a) = 0$,
with $a = (k, \varphi)$.

The idea behind our data generation is to randomly generate
the solution $U$ and the function $K$,
and use this information to compute the associated source term $\Phi$.
For this strategy to be efficient,
we need to have a rough idea of the family of solutions we wish to approximate.
This strategy's advantage,
in contrast to generating $\Phi$,
is that it avoids having to solve the PDE for data generation.
Data generation is summarized in \cref{alg:data_generation}.
It depends on the choice of a random generator,
whose definition is the objective of the remainder of this paragraph.

\begin{algorithm}[!ht]
    \caption{Data generation}
    \label{alg:data_generation}
    \begin{algorithmic}[1]

        \medskip

        \Require
        $p \in \mathbb{N}$,
        random generator for $U$ and $K$,
        $N_\text{data}$ number of data points
        \Ensure input and output datasets $\mathcal{X}$ and $\mathcal{Y}$

        \medskip

        \For {$i \in \{1, \dots, N_\text{data}\}$}
        \State randomly generate $U$ and $K$
        \State project $U$ and $K$
        onto the discrete space to obtain $u$ and $k$
        \State set $\varphi = E(u; k)$
        \State add ${X}_i=(k, \varphi)$
        to the input dataset $\mathcal{X}$ and
        ${Y}_i=u$ to the output dataset $\mathcal{Y}$
        \EndFor
    \end{algorithmic}
\end{algorithm}

To complete \cref{alg:data_generation},
a random generator to build interesting functions $U$ and $K$ has to be designed.
This problem depends on the chosen application.
Here, even though it is an important question,
we do not discuss how to construct a representative dataset for one given application.
Instead, we propose the following generic generator,
which serves as a proof of concept:
\begin{equation}
    \label{eq:generator_mixture}
    \mathtt{g}(x) =
    0.5 +
    \sum_{i=0}^n \mathcal{N}_{\mu_i,\sigma_i}(x),
\end{equation}
where $\mathcal{N}_{\mu,\sigma}$ is the probability density function
of the normal distribution, with mean $\mu$ and variance~$\sigma^2$.
The integer $n$ is randomly chosen between $0$ and $5$,
and thus the generated function~$\mathtt{g}$ corresponds to a sum of $n + 1$ Gaussian functions.
Each variance $\sigma_i^2$ is uniformly sampled in $[0.025,0.07]$,
and each mean $\mu_i$ is uniformly sampled in a ball of size $0.25$,
centered at the midpoint of the space domain $\Omega$.
This generator is then used to generate
both the solution $u$ and the space function $K$.
In principle, the method could work for other function families.
However, the larger the family, the poorer the performance.
The definition of this generator completes \cref{alg:data_generation}.
It is then used to define three datasets:
\begin{itemize}
    \item the training dataset,
          which includes several fixed mesh sizes;
    \item the validation dataset,
          used to compare the model to data and
          to select the hyperparameters,
          is constructed in the same way but using additional mesh sizes;
    \item the test dataset, used to evaluate the performance of the model,
          is constructed like the validation dataset.
\end{itemize}

\subsection{Neural network structure}
\label{sec:NN_structure}

\begin{figure}[ht!]
    \centering
    \begin{tikzpicture}

        \pgfmathsetmacro{\dx}{0.75}
        \pgfmathsetmacro{\dy}{-2}
        \pgfmathsetmacro{\df}{-4}

        \node[left, rectangle, draw] at (0,0) (input)
        {$\makecell{A(x) \\ \in \mathbb{R}^{n_c}}$};

        \draw [-latex] (input.east) -- ++(\dx,0)
        node(extrapolation) [right, rectangle, draw]
        {$\makecell{V(x) \\ \in \mathbb{R}^{n_p}}$}
        node [midway, below] {$P$};

        \draw [-latex] (extrapolation.east) -- ++(\dx,0)
        node(fourier_1) [right, rectangle, draw]
        {\makecell{Fourier \\ layer $1$}};

        \draw [-latex] (fourier_1.east) -- ++(\dx,0)
        node(fourier_2) [right, rectangle, draw]
        {\makecell{Fourier \\ layer $2$}};

        \draw [-latex] (fourier_2.east) -- ++(\dx,0)
        node(dots) [right]
        {\dots};

        \draw [-latex] (dots.east) -- ++(\dx,0)
        node(fourier_L) [right, rectangle, draw]
        {\makecell{Fourier \\ layer $L$}};

        \draw [-latex] (fourier_L.east) -- ++(\dx,0)
        node(projection) [right, rectangle, draw]
        {$\makecell{W(x) \\ \in \mathbb{R}^{n_p}}$};

        \draw [-latex] (projection.east) -- ++(\dx,0)
        node(output) [right, rectangle, draw]
        {$\makecell{U(x) \\ \in \mathbb{R}}$}
        node [midway, below] {$Q$};

        \draw [densely dotted]
        (fourier_2.south west) -- ($(input.west) + (0,\dy)$)
        coordinate[right] (corner_left);

        \draw [densely dotted]
        (fourier_2.south east) -- ($(output.east) + (0,\dy)$)
        coordinate[right] (corner_right);

        \draw
        (corner_left) -- (corner_right) --
        ($(corner_right) + (0, \df)$) --
        ($(corner_left) + (0, \df)$) --
        cycle node[midway, right, rectangle, draw, xshift=6pt] (fourier_input)
        {$\makecell{V(x) \\ \in \mathbb{R}^{n_p}}$};

        \node (fft_ifft) [rectangle, draw] at ($(fourier_input) + (5.5, 1)$) {
            \begin{tikzpicture}
                \node[left, rectangle, draw] at (0,0) (fft)
                {$\makecell{\mathbb{F}(V) \\ \in \mathbb{R}^{n_p/2}}$};

                \draw [-latex] (fft.east) -- ++(\dx,0)
                node(reduction) [right, rectangle, draw]
                {$\makecell{\widehat R(\mathbb{F}(V)) \\ \in \mathbb{R}^{m}}$};

                \draw [-latex] (reduction.east) -- ++(\dx,0)
                node(ifft) [right, rectangle, draw]
                {$\makecell{\mathbb{F}^{-1}(R(\mathbb{F}(V))) \\ \in \mathbb{R}^{n_p}}$};
            \end{tikzpicture}
        };

        \draw [-latex] (fourier_input.east) -- (fft_ifft.west);

        \node (linear) [rectangle, draw] at ($(fft_ifft) + (0, -2)$)
        {$\makecell{W(V(x)) \\ \in \mathbb{R}^{n_p}}$};

        \draw [-latex] (fourier_input.east) -- (linear.west);

        \node (plus) [rectangle, draw] at ($(fft_ifft) + (5.5, -1)$)
        {\sbox0{$\makecell{V(x) \\ \in \mathbb{R}^{n_p}}$}\makebox(\wd0,\ht0)[c]{$+$}};

        \draw [-latex] (fft_ifft.east) -- (plus.west);
        \draw [-latex] (linear.east) -- (plus.west);

        \node (sigma) [rectangle, draw] at ($(plus.east) + (1.5, 0)$)
        {\sbox0{$\makecell{V(x) \\ \in \mathbb{R}^{n_p}}$}\makebox(\wd0,\ht0)[c]{$\sigma$}};

        \draw [-latex] (plus.east) -- (sigma.west);

    \end{tikzpicture}
    \caption{%
        Sketch of an FNO, adapted from \cite{li2020fourier}.
        The network has $n_c$ input channels ($n_c = 2$ in 1D and $n_c = 5$ in 2D),
        which are then extrapolated to $n_p$ dimensions
        by the extrapolation layer $P$;
        $L$ Fourier layers follow,
        whose result is projected back to $\mathbb{R}$
        to give an approximation of the PDE solution.
        In the Fourier layers, $\mathbb{F}$ denotes
        the fast Fourier transform $\mathtt{FFT}$.
    }
    \label{fig:fno_archi}
\end{figure}
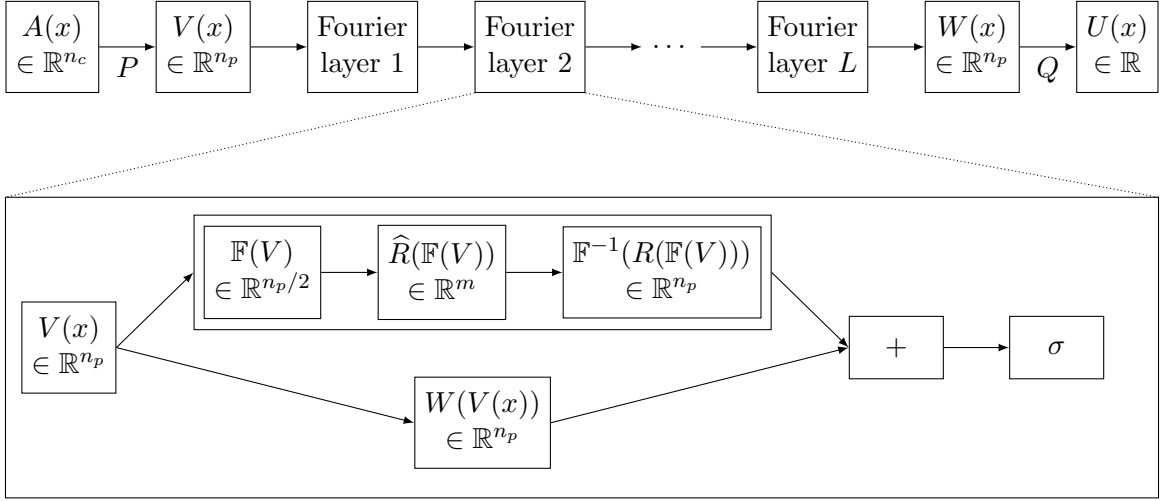

As mentioned before,
we elect to use an FNO in both cases \eqref{1deq} and \eqref{2deq}.
For the sake of completeness, we recall the main ideas behind FNOs,
and represent the architecture,
following \cite{li2020fourier}, on \cref{fig:fno_archi}.
The FNO is made of an extrapolation layer,
which starts by lifting the input to a higher dimension,
followed by several Fourier layers,
and ends with a projection layer,
projecting the output back to the original dimension.
Such a neural network is exactly equivalent to a fully convolutional one:
the input $X$ is transformed via a succession of convolutions
and non-linear activation functions,
but convolutions are performed via a Fast Fourier Transform algorithm:
\begin{equation*}
    \mathtt{FFT}^{-1}\Big( \widehat R \big(\mathtt{FFT}(X)\big) \Big),
\end{equation*}
where $\widehat R$ is a trainable vector with a small support $[0,m]$
(where e.g. $m=20$),
which multiplicatively modifies the low frequencies (modes) of the signal.
We can imagine that $\widehat R$ is itself the Fourier transform
of a (never computed) kernel $R$;
so our convolution is equivalent to the classical one $X\star R$.
In classical convolutional networks,
the kernels $R$ have very small support (e.g. 3),
while in the case of FNOs, the kernel has full spatial support.
For the full description of the model, we refer the reader to \cite{li2020fourier}.

The network always takes as inputs $\Phi$ and $K$;
in 1D, this corresponds to two channels,
and five channels in 2D (since $K$ is matrix-valued in this case).
The output is always the real-valued solution~$U$.
To obtain $G^+_\theta$,
this neural network is then batched on the $N_h$ discretization points,
to obtain a map from $\R^{n_c N_h}$ to $\R^{N_h}$.

This architecture involves several hyperparameters,
written on \cref{fig:fno_archi}:
the number $L$ of Fourier layers,
the number $m$ of Fourier modes to keep, and
the width $n_p$ (corresponding to the higher dimension after extrapolation).
To choose these hyperparameters, we perform a grid search, described below.

\subsection{Hyperparameter choice}
\label{sec:hyperparameter_choice}

So far, this section has been dedicated to the construction of
the neural network $G_\theta^+$ to provide a relevant initial guess,
and of the loss functions to be minimized during training.
The training process itself is rather standard,
based on gradient methods and back-propagation,
as is usual for neural networks.
However, this process depends on quite a lot of hyperparameters,
given in \cref{sec:hyperparameters};
it is crucial to choose them carefully to obtain a good performance.
To that end, we perform a grid search,
described in \cref{sec:grid_search}, to find the best hyperparameters.

\subsubsection{Hyperparameters}
\label{sec:hyperparameters}

As a first step, we recall all the hyperparameters involved in the learning process.
We split them into two categories:
the hyperparameters used for the training of the neural network,
and the hyperparameters of the neural network itself.
For each hyperparameter, we also give the possible values
explored in the grid search.

\begin{table}[ht!]
    \centering
    \begin{tabular}{lll}
        \toprule
        Training hyperparameters         &  & Set of values for the grid search \\
        \cmidrule(lr){1-3}
        $\ell_0$ (Initial learning rate) &  & $10^{-2},10^{-3},10^{-4}$         \\
        $\gamma$ (Exponential decay)     &  & $0.98,0.99,1$                     \\
        $N_b$ (Batch size)               &  & $64,128,256,512$                  \\
        $\omega$ (Loss function weight)  &  & $0,0.5,1$                         \\
        \bottomrule
    \end{tabular}
    \label{tab:hyp1}
    \caption{%
        Hyperparameters used for the training of the neural network,
        and values explored in the grid search.
    }
\end{table}

First, the training hyperparameters are collected in \cref{tab:hyp1}.
The first hyperparameter, $\ell_0$,
is the initial learning rate of the ADAM optimizer.
The learning rate is then exponentially decayed, with rate $\gamma$.
The parameter $N_b$ is the batch size,
i.e. the number of samples used to compute the loss functions.
The parameter $\omega$ corresponds to the weight between the loss functions.
The choice of the final loss function
will be discussed in more detail in \cref{sec:numerical_results}.
For the moment, suffice it to say that the loss function
combines the three loss functions
\eqref{eq_loss1}--\eqref{eq_loss2}--\eqref{eq_loss3}
differently in 1D, for \eqref{1deq}, and in 2D, for \eqref{2deq}.
In 1D, we use a combination between the $L^2$ data loss function \eqref{eq_loss1}
and the discretization-informed loss function \eqref{eq_loss3};
the resulting loss function is given by
\begin{equation}
    \label{eq:1D_loss}
    \mathcal{L}_\text{1D}(\theta) =\omega\mathcal{L}_{\text{data}}^{L^2}(\theta)+(1-\omega)10^{-4}\mathcal{L}_\text{dis}(\theta).
\end{equation}
In 2D, the combination is between the $L^2$ data loss function \eqref{eq_loss1}
and the $H^1$ data loss function~\eqref{eq_loss2},
to obtain the following loss function:
\begin{equation}
    \label{eq:2D_loss}
    \mathcal{L}_\text{2D}(\theta) =
    \omega\mathcal{L}_{\text{data}}^{L^2}(\theta)+
    (1-\omega)10^{-2}\mathcal{L}_{\text{data}}^{H^1}(\theta).
\end{equation}
The factors $10^{-4}$ in \eqref{eq:1D_loss} and $10^{-2}$ in \eqref{eq:2D_loss}
are introduced to normalize the discretization-informed and the $H^1$ loss functions,
which are always larger than the $L^2$ loss function.
Indeed, recall that the generator \eqref{eq:generator_mixture}
generates functions based on the probability density function
of a normal distribution.
Therefore, their spatial derivatives, and especially their second derivatives, have a norm greater than that of the functions themselves.
We have added the normalizing factors to account for this discrepancy.

\begin{table}[ht!]
    \centering
    \begin{tabular}{lll}
        \toprule
        Network hyperparameters &  & Set of values for the grid search \\
        \cmidrule(lr){1-3}
        $L$ (Fourier layers)    &  & $2,3,4,5$                         \\
        $m$ (Fourier modes)     &  & $10,20,30,40$                     \\
        $n_p$ (FNO width)       &  & $10,20,30,50$                     \\
        \bottomrule
    \end{tabular}
    \caption{%
        Hyperparameters used for the neural network itself,
        and values explored in the grid search.
    }
    \label{tab:hyp2}
\end{table}

Second, the hyperparameters of the neural network are summarized in \cref{tab:hyp2}.
We refer the reader to \cref{sec:NN_structure} and \cref{fig:fno_archi}
for an explanation of these hyperparameters.

\subsubsection{Grid search}
\label{sec:grid_search}

The main idea behind grid search is to run the training process
for all possible combinations of hyperparameters,
and to select the best combination.
We denote by $\bs{\mu}$ the set of hyperparameters of our problem:
\[
    \bs{\mu}= \{\ell_0,\gamma,N_b,\omega,L,m,p\}.
\]
Grid search can be formalized as the following optimization problem:
\[
    \bs{\mu}_\text{opt} = \argmin_{\bs{\mu}} S(\bs{\mu}),
\]
with $S$ a score function to be defined and
$\bs{\mu}_\text{opt}$ the optimal parameters.
To avoid overfitting and increase the generalization ability of our model,
this score function will be computed on validation data,
which is generated following \cref{sec:data_generation}.
We denote by $N_\text{val}$ the number of validation data.
As we will see upon defining a suitable score function,
this score function will not necessarily be
differentiable with respect to $\bs{\mu}$.
This is why we choose to use a grid search rather than
a gradient-based optimization algorithm.

The classical choice in the machine learning community is to
set the score function equal to parts of the loss function of the problem,
given here by \eqref{eq:1D_loss} or \eqref{eq:2D_loss}.
We thus define two score functions:
$S_\text{data}$, based on the $L^2$ loss function,
and $S_\text{dis}$, based on the discretization-informed loss function.
Their evaluations are summarized in \cref{alg:score_function_with_loss}.

\begin{algorithm}[!ht]
    \caption{Evaluation of $S_\text{data}(\bs{\mu})$ and $S_\text{dis}(\bs{\mu})$}
    \label{alg:score_function_with_loss}
    \begin{algorithmic}[1]

        \medskip

        \Require
        set of parameters $\bs{\mu}$,
        $n_\text{epoch}$ number of training epochs,
        $N_\text{val}$ number of validation data,
        $(X_i, Y_i)_{i \in \{1, \dots, N_\text{val}\}}$ validation data
        \Ensure
        scores $S_\text{data}(\bs{\mu})$ and $S_\text{dis}(\bs{\mu})$

        \medskip

        \State train $G_\theta^+$ during $n_\text{epoch}$ epochs
        with hyperparameters $\bs{\mu}$
        \State compute the scores
        $\displaystyle
            S_\text{data}(\bs{\mu})
            =
            \sum_{i=1}^{N_\text{val}}\| G_\theta^+({X}_i)-{Y}_i \|_2^2
        $ and
        $\displaystyle
            S_\text{dis}(\bs{\mu})
            =
            \sum_{i=1}^{N_\text{val}}\| F(G_\theta^+({X}_i){\rb;{X}_i)} \|_2^2
        $
    \end{algorithmic}
\end{algorithm}

However, using such a score function implies the underlying assumption that
the better the neural network approximation of the solution,
the faster the convergence of Newton's method.
Even though this assumption seems natural, it remains debatable.
Indeed, there could be cases where the model provides a good approximation
of the PDE solution in the $L^2$ norm,
but local errors or irregularities could slow down the convergence of Newton's method
if using this approximation as an initial guess.
Since our end goal is to accelerate the convergence of Newton's method,
another natural choice for the score function
would be to base it on the actual number of iterations.
Computing such a score function amounts to
running Newton's algorithm with the neural network prediction
and with the naive initial guess, and comparing the number of iterations.
The score is then defined as
the average gain in iterations between the two approaches.
As a consequence, we suggest computing this new score function
$S_\text{iter}(\bs{\mu})$ via \cref{alg:score_function_based_on_iterations}.
{\ra%
Note that, in this algorithm,
we set a maximum number $M_\text{iter}$ of Newton iterations
to avoid the algorithm running indefinitely.
In practice, we set $M_\text{iter} = 2000$.
}

\begin{algorithm}[!ht]
    \caption{Evaluation of $S_\text{iter}(\bs{\mu})$}
    \label{alg:score_function_based_on_iterations}
    \begin{algorithmic}[1]

        \medskip

        \Require
        set of hyperparameters $\bs{\mu}$,
        $n_\text{epoch}$ number of training epochs,
        $N_\text{val}$ number of validation data,
        $(X_i)_{i \in \{1, \dots, N_\text{val}\}}$ validation data,
        {\ra $M_\text{iter}$ maximum number of Newton iterations}
        \Ensure
        score $S_\text{iter}(\bs{\mu})$

        \medskip

        \State train $G_\theta^+$ during $n_\text{epoch}$ epochs
        with hyperparameters $\bs{\mu}$
        \For {$i \in \{1, \dots, N_\text{val}\}$}
        \State compute the initial guess $u_0=G_\theta^+({X}_i)$
        using the trained neural network
        \State solve the problem with a JFNK solver,
        with initial guess $u_0$,
        and denote by $k_i$ the number of iterations
        needed to reach convergence{\ra; if $k_i > M_\text{iter}$, set $k_i = M_\text{iter}$}
        \State solve the problem with a JFNK solver,
        with the naive initial guess $(1, 1, \dots, 1)$,
        and denote by $\tilde k_i$ the number of iterations
        needed to reach convergence{\ra; if $\tilde k_i > M_\text{iter}$, set $k_i = M_\text{iter}$}
        \EndFor
        \State compute the score
        $\displaystyle S_\text{iter}(\bs{\mu})=\frac{1}{N_\text{val}}\sum_{i=1}^{N_\text{val}} \frac{\tilde{k}_i}{k_i}$
    \end{algorithmic}
\end{algorithm}

To evaluate the best values of each hyperparameter,
we evaluated the three scores on several mesh resolutions.
For the loss-based scores $S_\text{data}$ and $S_\text{dis}$,
we used $N_\text{val}=1000$ validation data for each mesh resolution.
The iteration-based score $S_\text{iter}$, however,
is more expensive to compute since Newton's method
has to be run up to convergence.
Therefore, to compute this score, performed $N_\text{val}=5$ Newton resolutions
per mesh resolution.

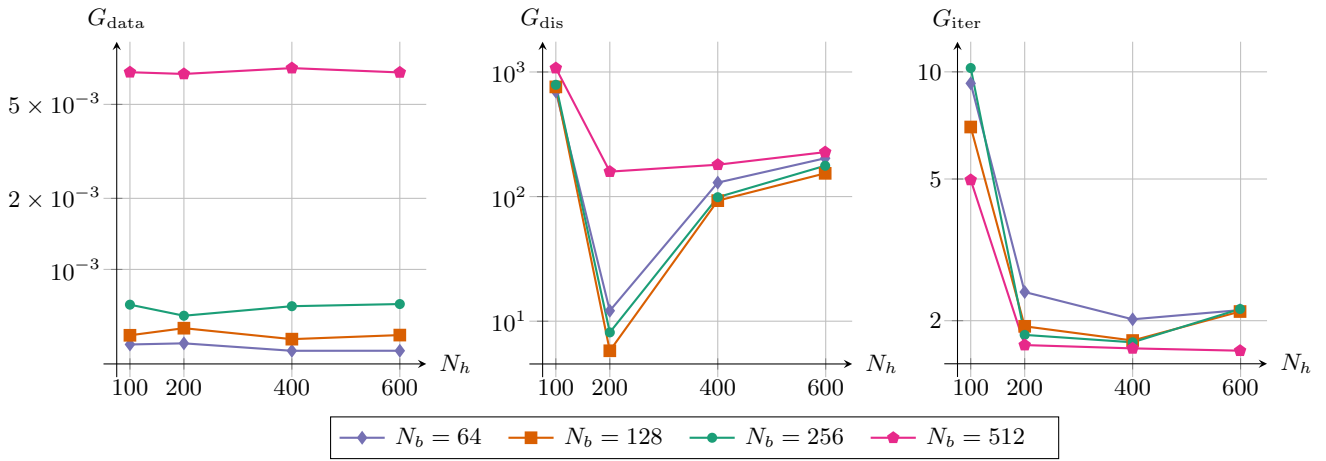
\begin{figure}[ht!]
    \centering
    \setlength{\plotwidth}{0.34\linewidth}
    \setlength{\plotheight}{0.35\linewidth}

    \pgfplotsset{
        errors_wrt_batch_size/.style={
                axis lines = left,
                enlarge x limits={abs=10pt},
                enlarge y limits={abs=10pt},
                axis x line shift = -5pt,
                axis y line shift = -5pt,
                xlabel = {$N_h$},
                xlabel style={at={(ticklabel* cs:1.01)},anchor=west,font=\footnotesize,},
                ylabel style={at={(ticklabel* cs:1.01)},anchor=west,font=\footnotesize,},
                xmin = 100, xmax = 600,
                xtick={100, 200, 400, 600},
                tick label style={font=\footnotesize},
                yminorticks=false,
                width=\plotwidth, height=\plotheight,
                mark options={solid, scale=1},
                grid = major,
                legend style={font=\footnotesize},
            }
    }

    \begin{tikzpicture}

        \begin{semilogyaxis}[
                errors_wrt_batch_size,
                name = left_plot,
                ylabel = {\rotatebox{270}{$G_\text{data}$}},
                ytick={1e-3, 2e-3, 5e-3},
                yticklabels={$10^{-3}$, $2 \times 10^{-3}$, $5 \times 10^{-3}$},
                legend entries={\,$N_b = 64$\;\;\;, \,$N_b = 128$\;\;\;, \,$N_b = 256$\;\;\;, \,$N_b = 512$\;\;\;},
                legend columns=4,
                legend to name=leg:errors_wrt_batch_size,
                legend image post style={mark options={solid, scale=1}},
            ]

            \addplot [style={solid}, mark=diamond*,  mark size=2, line width=0.8pt, color=graph_1 ]
            table [x=nx, y=64, col sep=comma]
                {plots/data_error_batchsize.csv};

            \addplot [style={solid}, mark=square*,   mark size=2, line width=0.8pt, color=graph_2 ]
            table [x=nx, y=128, col sep=comma]
                {plots/data_error_batchsize.csv};

            \addplot [style={solid}, mark=*, mark size=1.5, line width=0.8pt, color=graph_3]
            table [x=nx, y=256, col sep=comma]
                {plots/data_error_batchsize.csv};

            \addplot [style={solid}, mark=pentagon*, mark size=2, line width=0.8pt, color=graph_4 ]
            table [x=nx, y=512, col sep=comma]
                {plots/data_error_batchsize.csv};

        \end{semilogyaxis}

        \begin{semilogyaxis}[
                errors_wrt_batch_size,
                name = middle_plot,
                at=(left_plot.outer east), anchor=outer west,
                ylabel = {\rotatebox{270}{$G_\text{dis}$}},
            ]

            \addplot [style={solid}, mark=diamond*,  mark size=2, line width=0.8pt, color=graph_1 ]
            table [x=nx, y=64, col sep=comma]
                {plots/residu_compared_batchsize.csv};

            \addplot [style={solid}, mark=square*,   mark size=2, line width=0.8pt, color=graph_2 ]
            table [x=nx, y=128, col sep=comma]
                {plots/residu_compared_batchsize.csv};

            \addplot [style={solid}, mark=*, mark size=1.5, line width=0.8pt, color=graph_3]
            table [x=nx, y=256, col sep=comma]
                {plots/residu_compared_batchsize.csv};

            \addplot [style={solid}, mark=pentagon*, mark size=2, line width=0.8pt, color=graph_4 ]
            table [x=nx, y=512, col sep=comma]
                {plots/residu_compared_batchsize.csv};

        \end{semilogyaxis}

        \begin{semilogyaxis}[
                errors_wrt_batch_size,
                name = right_plot,
                at=(middle_plot.outer east), anchor=outer west,
                ylabel = {\rotatebox{270}{$G_\text{iter}$}},
                ytick={1,2,5,10},
                yticklabels={$1$, $2$, $5$, $10$},
            ]

            \addplot [style={solid}, mark=diamond*,  mark size=2, line width=0.8pt, color=graph_1 ]
            table [x=nx, y=64, col sep=comma]
                {plots/time_ratio_compared_batchsize.csv};

            \addplot [style={solid}, mark=square*,   mark size=2, line width=0.8pt, color=graph_2 ]
            table [x=nx, y=128, col sep=comma]
                {plots/time_ratio_compared_batchsize.csv};

            \addplot [style={solid}, mark=*, mark size=1.5, line width=0.8pt, color=graph_3]
            table [x=nx, y=256, col sep=comma]
                {plots/time_ratio_compared_batchsize.csv};

            \addplot [style={solid}, mark=pentagon*, mark size=2, line width=0.8pt, color=graph_4 ]
            table [x=nx, y=512, col sep=comma]
                {plots/time_ratio_compared_batchsize.csv};

        \end{semilogyaxis}

        \node[yshift=-12.5pt] at (middle_plot.outer south) {\pgfplotslegendfromname{leg:errors_wrt_batch_size}};

    \end{tikzpicture}
    \caption{%
        Score with respect to the number of points in the grid,
        for different values of the batch size~$N_b$.
        From left to right, we display scores
        $S_\text{data}$, $S_\text{dis}$ and $S_\text{iter}$.
        Scores $S_\text{data}$ and $S_\text{dis}$
        should be as small as possible;
        score $S_\text{iter}$ should be as large as possible.
    }
    \label{fig:compare_bT}
\end{figure}

On \cref{fig:compare_bT}, we display the results of the grid search
on the 1D problem \eqref{1deq},
only for the batch size hyperparameter $N_b$ for simplicity.
From left to right,
we display scores $S_\text{data}$, $S_\text{dis}$ and $S_\text{iter}$,
with respect to the number $N_h$ of points,
for different values of the batch size $N_b$.
The left and center panels show that a batch size of $N_b=512$ is not efficient,
while the other batch sizes yield comparable scores.
Moving on to the right panel, we observe that $N_b=64$ is the best choice
to reduce the number of iterations of Newton's method.
Consequently, from these indicators, we choose $N_b=64$.

    {\rb%
        All in all, the hyperparameters are chosen by performing the grid search.
        Note that, in most cases,
        all three scores agree on the best value of each hyperparameter.
        In case the three scores did not concur,
        we have used the hyperparameters associated with $S_\text{iter}$.
        Indeed, this score is the one deemed to be the most important
        in practice,
        since it is associated with reducing the number of iterations.%
    }
The final choice of all hyperparameters is summarized in \cref{tab:hyp_final}.

\begin{table}[ht!]
    \centering
    \begin{tabular}{llc}
        \toprule
        Hyperparameters                  &  & Chosen value \\
        \cmidrule(lr){1-3}
        $\ell_0$ (Initial learning rate) &  & $10^{-3}$    \\
        $\gamma$ (Exponential decay)     &  & $0.99$       \\
        $N_b$ (Batch size)               &  & $64$         \\
        $\omega$ (Loss function weight)  &  & $0.5$        \\
        $L$ (Fourier layers)             &  & $4$          \\
        $m$ (Fourier modes)              &  & $30$         \\
        $n_p$ (FNO width)                &  & $30$         \\
        \bottomrule
    \end{tabular}
    \caption{%
        Summary of the hyperparameters;
        for each hyperparameter, we give the best value obtained through grid search.
    }
    \label{tab:hyp_final}
\end{table}

\section{Numerical results}
\label{sec:numerical_results}

Equipped with the values of the hyperparameters
collected in \cref{tab:hyp_final},
we now present some numerical results to validate our approach.
Namely, we compare the convergence of Newton's method
using a naive, constant initial guess
and using our FNO as a predicted initial guess.
We first tackle the 1D case in \cref{sec:1D_results},
and move on to the 2D case in \cref{sec:2D_results}.
In each case, we start by displaying the results of the FNO prediction,
to make sure they are close to the exact solution.
Then, we compare both types of initial guesses.

Throughout this section,
we use the \texttt{newton\_krylov} solver from the \texttt{scipy} library.
It is already a very efficient solver, with built-in advanced optimizations.
It is a Jacobian-free Newton-Krylov method,
where linear solves are performed using the LGMRES method
and where the convergence is enhanced with
an Armijo line search method {\ra from \cite{PhilipCONVERGENCECF}}.
Such a method is typical for large-scale problems,
for instance, ones stemming from PDE discretizations.
    {\ra%
        Moreover, since this work is essentially a proof of concept,
        we remain on relatively coarse meshes in 2D (up to $100^2$ points),
        while finer meshes are considered in 1D (up to $600$ points).
    }

\subsection{Results in one space dimension}
\label{sec:1D_results}

The first batch of experiments we run are based on the 1D PDE \eqref{1deq},
with the nonlinearity $p$ fixed to $4$.
The coefficient $\alpha_0$ tunes the amplitude of the nonlinearity, larger values of $\alpha_0$ lead to a prevalence of the nonlinear part, and thus make the JFNK method slower to converge.
We discretize this equation with a classical finite difference scheme.

For each experiment, we fix a value of $\alpha_0$ to represent the difficulty of the problem.
The model is trained on two meshes made of $200$ and $400$ points.
These experiments aim at checking the ability of the FNO to predict a solution for an unknown mesh resolution.
    {\rall We first motivate, in \cref{sec:motivation},
        the need for suitable initial guesses,
        even in the case of 1D equations such as \eqref{1deq}.}
We then discuss the choice of the loss function in \cref{sec:1D_prediction},
and check the relevance of the initial guess provided by the FNO
for a small nonlinearity $\alpha_0=2$ in \cref{sec:1D_small_alpha}
and for larger nonlinearities $\alpha_0 \in \{5, 8\}$ in \cref{sec:1D_large_alpha}.
{\ra \cref{sec:1D_large_alpha} also contains
a study of the number of simulations needed
to overcome the training time of the FNO.}
{\rall We finally compare our approach to other methods
from the literature in \cref{sec:comparison}.}

{\rall
\subsubsection{Motivation: failure of classical initial conditions for Newton's method}
\label{sec:motivation}

First, we illustrate the need for a suitable initial guess
to apply Newton's method to solve the discrete
1D PDE \eqref{1deq}.
We consider meshes made of
$N_h = 100$ points in \cref{tab:initializations_100_pts},
$N_h = 200$ points in \cref{tab:initializations_200_pts},
$N_h = 400$ points in \cref{tab:initializations_400_pts},
and $N_h = 600$ points in \cref{tab:initializations_600_pts}.
For each mesh, we take values of $\alpha_0$ in $\{2, 5, 8\}$,
ranging from weak to strong nonlinearities.
In such elliptic problems,
we usually do not have any prior on the solution,
so changing initial guesses usually amounts
to changing the constant value used as an initial guess,
potentially adding some random noise.
We consider several initial conditions for the Newton's method:
\begin{itemize}[nosep]
    \item five constant initial conditions
          (with values $0.5$, $0.75$, $1$, $1.25$, $1.5$);
    \item a constant initialization with value $1$
          added to a Gaussian perturbation $\mathcal{N}_{0, 0.25}$,
          with mean $0$ and standard deviation $0.25$;
    \item a constant initialization with value $1$
          added to a Gaussian perturbation
          $\mathcal{N}_{0,  \sigma_\text{exact}}$,
          where $\sigma_\text{exact}$ is the variance
          of the exact solution;
    \item and the exact solution $U_\text{exact}$
          added to a Gaussian perturbation $\mathcal{N}_{0, 0.25}$.
\end{itemize}
Note that the last two initial conditions are not realistic,
since they depend on prior knowledge of the exact solution.

\begin{table}[!ht]
    \rall
    \centering
    \begin{tabular}{ccccccc}
        \toprule
        \multirow{2}{*}[-2pt]{\makecell{initial                                                                                                                            \\ condition}}
                                                   & \multicolumn{2}{c}{$\alpha_0=2$}
                                                   & \multicolumn{2}{c}{$\alpha_0=5$}
                                                   & \multicolumn{2}{c}{$\alpha_0=8$}                                                                                      \\
        \cmidrule{2-7}
                                                   & $\bar N_\text{iter}$             & \% failure & $\bar N_\text{iter}$ & \% failure & $\bar N_\text{iter}$ & \% failure \\
        \cmidrule(lr){1-1}\cmidrule(lr){2-3}\cmidrule(lr){4-5}\cmidrule(lr){6-7}
        $0.5$                                      & 2000                             & 100        & 2000                 & 100        & 2000                 & 100        \\
        $0.75$                                     & 1686                             & 72.5       & 1670                 & 72.5       & 1804                 & 80         \\
        $1$                                        & 280                              & 0          & 1293                 & 30         & 1824                 & 75         \\
        $1.25$                                     & 1435                             & 42.5       & 1993                 & 97.5       & 2000                 & 100        \\
        $1.5$                                      & 1960                             & 95         & 2000                 & 100        & 2000                 & 100        \\
        $1 + \mathcal{N}_{0, 0.25}$                & 362                              & 5          & 1244                 & 27.5       & 1791                 & 70         \\
        $1 + \mathcal{N}_{0, \sigma_\text{exact}}$ & 285                              & 0          & 1277                 & 32.5       & 1801                 & 75         \\
        $U_\text{exact} + \mathcal{N}_{0, 0.25}$   & 285                              & 10         & 338                  & 7.5        & 580                  & 20         \\
        \bottomrule
    \end{tabular}
    \caption{Comparison of initial conditions for different values of $\alpha_0$ and $N_h = 100$ points.}
    \label{tab:initializations_100_pts}
\end{table}

\begin{table}[!ht]
    \rall
    \centering
    \begin{tabular}{ccccccc}
        \toprule
        \multirow{2}{*}[-2pt]{\makecell{initial                                                                                                                            \\ condition}}
                                                   & \multicolumn{2}{c}{$\alpha_0=2$}
                                                   & \multicolumn{2}{c}{$\alpha_0=5$}
                                                   & \multicolumn{2}{c}{$\alpha_0=8$}                                                                                      \\
        \cmidrule{2-7}
                                                   & $\bar N_\text{iter}$             & \% failure & $\bar N_\text{iter}$ & \% failure & $\bar N_\text{iter}$ & \% failure \\
        \cmidrule(lr){1-1}\cmidrule(lr){2-3}\cmidrule(lr){4-5}\cmidrule(lr){6-7}
        $0.5$                                      & 2000                             & 100        & 2000                 & 100        & 2000                 & 100        \\
        $0.75$                                     & 1923                             & 95         & 1902                 & 92.5       & 1951                 & 95         \\
        $1$                                        & 351                              & 0          & 781                  & 5          & 1288                 & 32.5       \\
        $1.25$                                     & 1184                             & 25         & 1943                 & 95         & 1990                 & 97.5       \\
        $1.5$                                      & 1955                             & 95         & 2000                 & 100        & 2000                 & 100        \\
        $1 + \mathcal{N}_{0, 0.25}$                & 394                              & 0          & 889                  & 2.5        & 1453                 & 40         \\
        $1 + \mathcal{N}_{0, \sigma_\text{exact}}$ & 356                              & 0          & 791                  & 0          & 1384                 & 32.5       \\
        $U_\text{exact} + \mathcal{N}_{0, 0.25}$   & 312                              & 5          & 304                  & 0          & 525                  & 7.5        \\
        \bottomrule
    \end{tabular}
    \caption{Comparison of initial conditions for different values of $\alpha_0$ and $N_h = 200$ points.}
    \label{tab:initializations_200_pts}
\end{table}

\begin{table}[!ht]
    \rall
    \centering
    \begin{tabular}{ccccccc}
        \toprule
        \multirow{2}{*}[-2pt]{\makecell{initial                                                                                                                            \\ condition}}
                                                   & \multicolumn{2}{c}{$\alpha_0=2$}
                                                   & \multicolumn{2}{c}{$\alpha_0=5$}
                                                   & \multicolumn{2}{c}{$\alpha_0=8$}                                                                                      \\
        \cmidrule{2-7}
                                                   & $\bar N_\text{iter}$             & \% failure & $\bar N_\text{iter}$ & \% failure & $\bar N_\text{iter}$ & \% failure \\
        \cmidrule(lr){1-1}\cmidrule(lr){2-3}\cmidrule(lr){4-5}\cmidrule(lr){6-7}
        $0.5$                                      & 2000                             & 100        & 2000                 & 100        & 2000                 & 100        \\
        $0.75$                                     & 2000                             & 100        & 2000                 & 100        & 2000                 & 100        \\
        $1$                                        & 1157                             & 5          & 1888                 & 70         & 2000                 & 100        \\
        $1.25$                                     & 1246                             & 2.5        & 1941                 & 85         & 2000                 & 100        \\
        $1.5$                                      & 1483                             & 27.5       & 1974                 & 92.5       & 2000                 & 100        \\
        $1 + \mathcal{N}_{0, 0.25}$                & 1307                             & 15         & 1959                 & 77.5       & 1997                 & 97.5       \\
        $1 + \mathcal{N}_{0, \sigma_\text{exact}}$ & 1181                             & 2.5        & 1931                 & 80         & 2000                 & 100        \\
        $U_\text{exact} + \mathcal{N}_{0, 0.25}$   & 642                              & 2.5        & 1033                 & 10         & 1249                 & 7.5        \\
        \bottomrule
    \end{tabular}
    \caption{Comparison of initial conditions for different values of $\alpha_0$ and $N_h = 400$ points.}
    \label{tab:initializations_400_pts}
\end{table}

\begin{table}[!ht]
    \rall
    \centering
    \begin{tabular}{ccccccc}
        \toprule
        \multirow{2}{*}[-2pt]{\makecell{initial                                                                                                                            \\ condition}}
                                                   & \multicolumn{2}{c}{$\alpha_0=2$}
                                                   & \multicolumn{2}{c}{$\alpha_0=5$}
                                                   & \multicolumn{2}{c}{$\alpha_0=8$}                                                                                      \\
        \cmidrule{2-7}
                                                   & $\bar N_\text{iter}$             & \% failure & $\bar N_\text{iter}$ & \% failure & $\bar N_\text{iter}$ & \% failure \\
        \cmidrule(lr){1-1}\cmidrule(lr){2-3}\cmidrule(lr){4-5}\cmidrule(lr){6-7}
        $0.5$                                      & 2000                             & 100        & 2000                 & 100        & 2000                 & 100        \\
        $0.75$                                     & 1987                             & 97.5       & 2000                 & 100        & 2000                 & 100        \\
        $1$                                        & 1736                             & 45         & 2000                 & 100        & 2000                 & 100        \\
        $1.25$                                     & 1934                             & 87.5       & 2000                 & 100        & 2000                 & 100        \\
        $1.5$                                      & 1982                             & 87.5       & 2000                 & 100        & 2000                 & 100        \\
        $1 + \mathcal{N}_{0, 0.25}$                & 1724                             & 42.5       & 2000                 & 100        & 2000                 & 100        \\
        $1 + \mathcal{N}_{0, \sigma_\text{exact}}$ & 1684                             & 35         & 2000                 & 100        & 2000                 & 100        \\
        $U_\text{exact} + \mathcal{N}_{0, 0.25}$   & 1018                             & 2.5        & 1463                 & 5          & 1739                 & 25         \\
        \bottomrule
    \end{tabular}
    \caption{Comparison of initial conditions for different values of $\alpha_0$ and $N_h = 600$ points.}
    \label{tab:initializations_600_pts}
\end{table}

In \cref{tab:initializations_100_pts,tab:initializations_200_pts,tab:initializations_400_pts,tab:initializations_600_pts},
we report the average number of iterations $\bar N_\text{iter}$
that Newton's method took to converge,
as well as the percentage of failures,
over $40$ different examples generated with the
strategy from \cref{alg:data_generation}.
In the case of a failure,
we set the number of iterations to $2000$,
which is also the maximum number of Newton iterations
$M_\text{iter}$.
We observe that, as the number $N_h$ of points
and the strength $\alpha_0$ of the nonlinearity grow,
each initialization method becomes less efficient.
For instance, in \cref{tab:initializations_600_pts},
all naive initial conditions have failed,
and much information needs to be added to the initial guess
to notably decrease the failure rate.
This is a clear indication that the classical initial conditions
of Newton's method are not suitable for this problem.
To remedy this issue, we show in the remainder of this section
that our approach can provide a better initial guess
and converge in all cases.
}

\subsubsection{Choice of the loss function}
\label{sec:1D_prediction}

In this section, we discuss the choice of the loss function
by checking the ability of the FNO to approximate
the solution to \eqref{1deq},
and to its discretized version \eqref{eq:discretized_pde}.
Since the goal is to compare the loss functions,
we train a smaller FNO, with $L=4$, $m=20$ and $p=20$.
Moreover, we take $\alpha_0=5$ to have a strong nonlinearity.

\begin{figure}[ht!]
    \centering
    \begin{tikzpicture}
        \node (fig) at (0,0) {\includegraphics[width=0.98\textwidth]{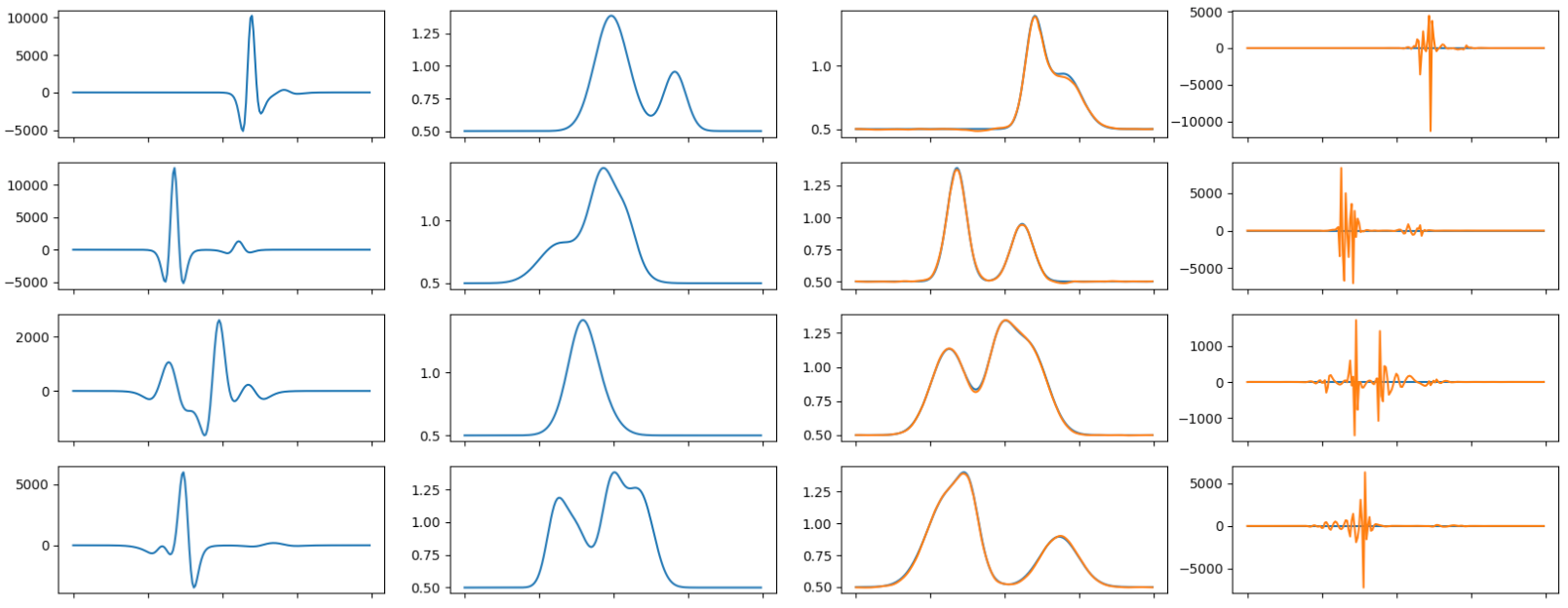}};
        \node [yshift=2pt] at (fig.151) {$\Phi(x)$};
        \node [yshift=2pt] at (fig.117) {$K(x)$};
        \node [yshift=2pt] at (fig.54) {$U(x)$};
        \node [yshift=2pt] at (fig.28) {Residual};
    \end{tikzpicture}
    \caption{%
        Examples obtained after training the FNO with the loss function $\mathcal{L}_\text{data}^{L^2}$ from \eqref{eq_loss2} only.
        From top to bottom, each series of graphs corresponds to a different example.
        In the first three columns, from left to right, we display $\Phi(x)$, $K(x)$,
        the solution $U(x)$ (blue line) and its prediction (orange line).
        The rightmost column compares the residuals of the discretized PDE computed with the exact solution (blue line) and the predicted solution (orange line).
    }
    \label{fig:1Dprediction}
\end{figure}
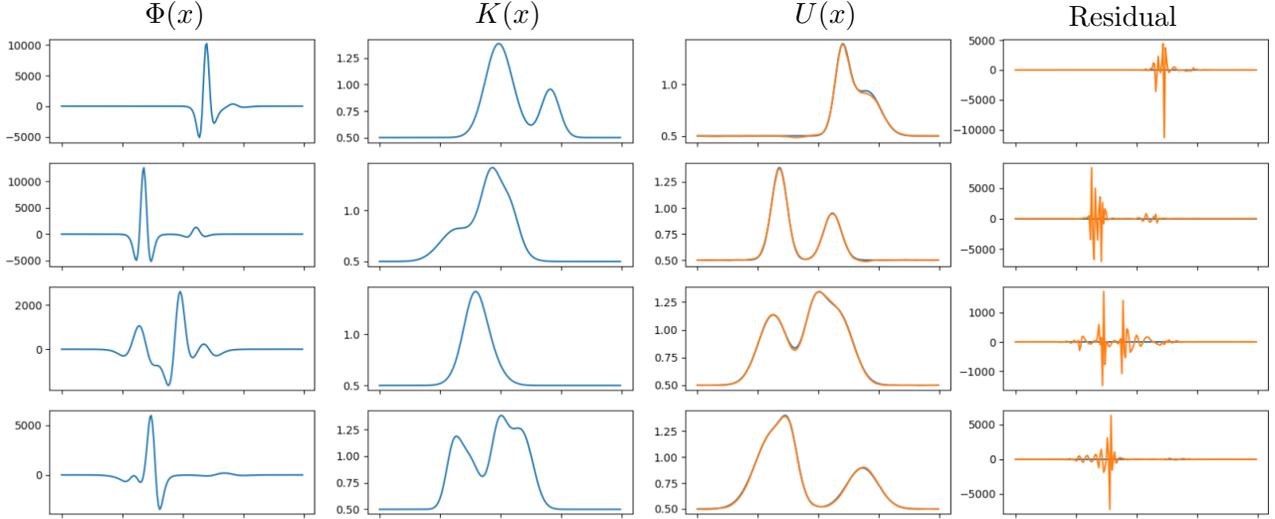

\begin{figure}[ht!]
    \centering
    \begin{tikzpicture}
        \node (fig) at (0,0) {\includegraphics[width=0.98\textwidth]{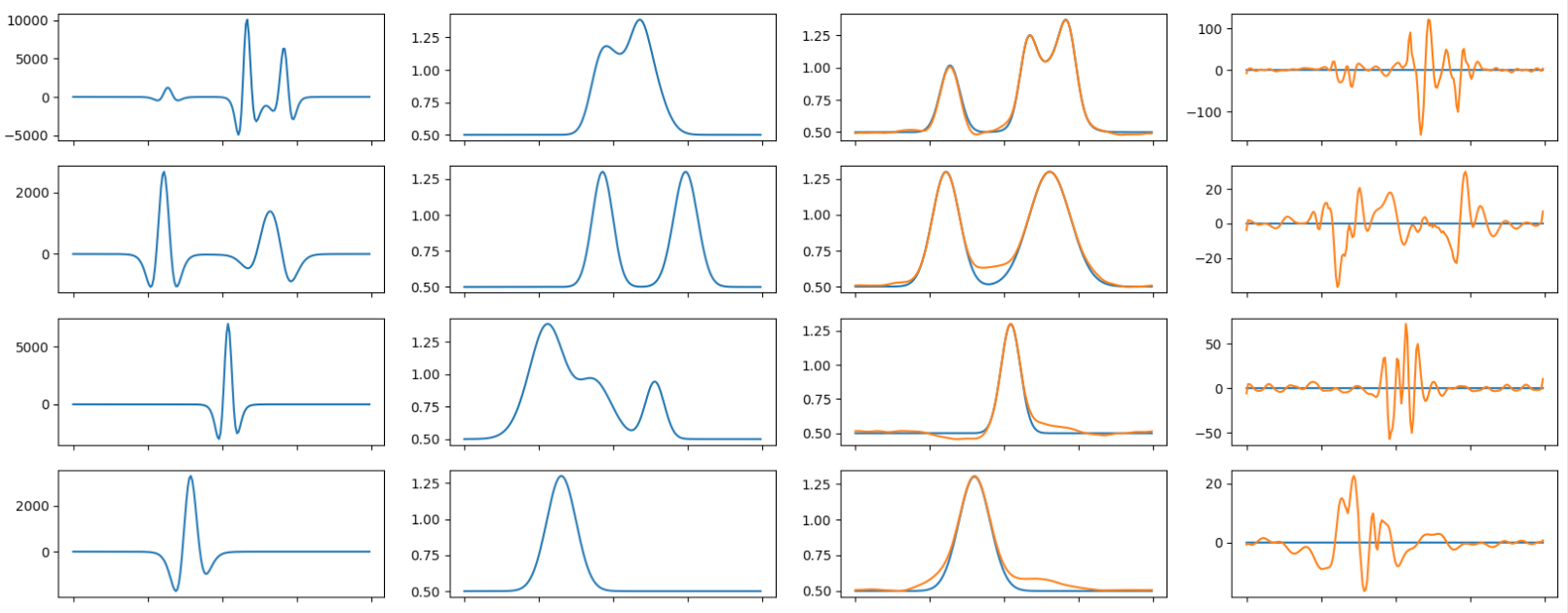}};
        \node [yshift=2pt] at (fig.151) {$\Phi(x)$};
        \node [yshift=2pt] at (fig.117) {$K(x)$};
        \node [yshift=2pt] at (fig.54) {$U(x)$};
        \node [yshift=2pt] at (fig.28) {Residual};
    \end{tikzpicture}
    \caption{%
        Examples obtained after training the FNO with
        the loss function $\mathcal{L}_\text{data}^{L^2}$ from \eqref{eq_loss2}
        and the discretization-informed loss
        function~$\mathcal{L}_\text{dis}$ from \eqref{eq_loss1}.
        From top to bottom,
        the four graphs correspond to a different example.
        In the first three columns,
        from left to right, we display $\Phi(x)$,~$K(x)$,
        the solution $U(x)$ (blue line) and its prediction (orange line).
        The rightmost column compares the residuals of the discretized PDE
        computed with the exact solution (blue line)
        and the predicted solution (orange line).
    }
    \label{fig:1Dprediction2}
\end{figure}
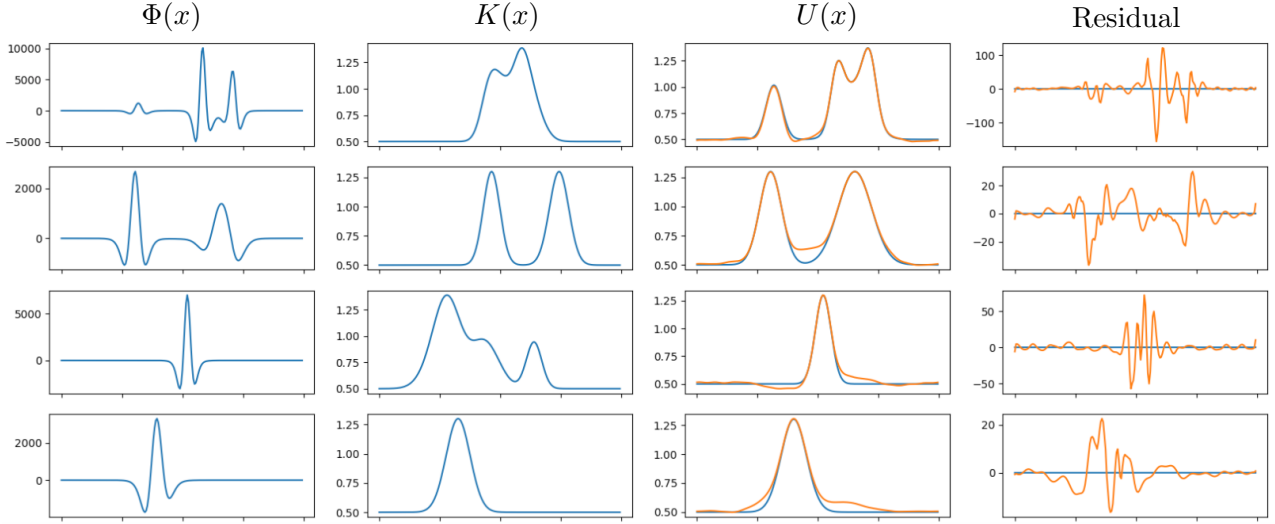

We compare two training strategies:
the first one uses only the $L^2$ data loss function given by \eqref{eq_loss2},
while the second one uses both the $L^2$ data loss function \eqref{eq_loss2}
and the discretization-informed loss function from~\eqref{eq_loss1}.
In both cases, the networks are trained
until the mean squared error (MSE) reaches around $2 \times 10^{-4}$.
These two strategies are compared in \cref{fig:1Dprediction,fig:1Dprediction2},
where $4$ random examples are displayed.
In each example, the functions $\Phi$ and $K$ are generated following \cref{sec:data_generation}, and the FNO is used to predict the resulting solution $u$.
In each row, we display the functions~$\Phi$,~$K$,~$u$ (predicted in blue and exact in orange) from left to right.
The rightmost figure compares the residuals of the discretized PDE, in the sense of \eqref{eq:general_pde}, computed with the exact solution (blue line) and the predicted solution (orange line).

In the case of \cref{fig:1Dprediction}, only the $L^2$ data loss function is used;
training time is below $10$ minutes on a cloud architecture with shared GPUs.
We note that the FNO provides a quite satisfactory prediction of the solution $U$.
However, we also remark that the residual obtained by plugging this solution
into the numerical scheme becomes quite large.
In light of the PDE \eqref{1deq},
since the third column shows that $U$ is correctly approximated by the FNO,
the high amplitude of the residual in the fourth column
can only be due to a poor approximation of
the discrete derivatives of~$U$ (and especially of its second derivative).

\cref{fig:1Dprediction2} shows the results obtained for an FNO  trained with both the $L^2$ data and discretization-informed loss function;
training time is below $15$ minutes,
still on a cloud architecture with shared GPUs.
On the one hand,
we observe that the network yields a good approximation of the solution $U$,
but that it performs noticeably worse than in the previous case in terms of approximation of $u$.
On the other hand, the residual of the scheme is much smaller,
which will be beneficial for the convergence of Newton's method.
In fact, by minimizing the discretization residual,
we have obtained a better reconstruction of the residual,
and therefore of the discrete derivatives of~$U$.
In practice, we will use this second training strategy for the final results,
because it produces a better initial guess for the JFNK method.

\subsubsection{\texorpdfstring{Results for $\alpha_0=2$}{Results for alpha0=2}}
\label{sec:1D_small_alpha}

We start by testing the approach with a weak nonlinearity ($\alpha_0=2$).
In this case, Newton's method is able to converge quite quickly
with the naive initial guess.
Nevertheless, we expect our predicted initial guess to outperform
the naive one by a significant margin.

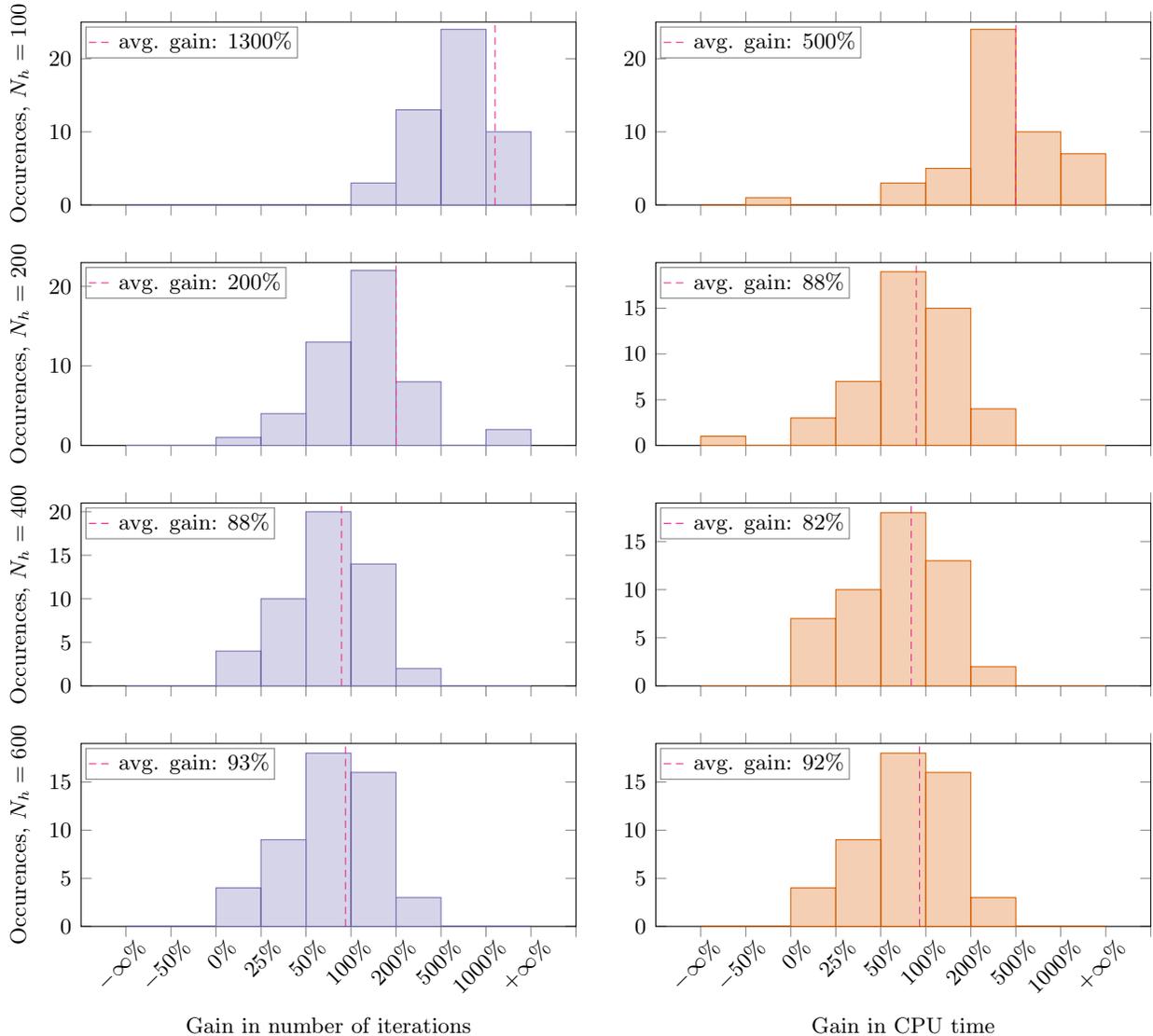
\begin{figure}[ht!]
    \centering
    \setlength{\plotwidth}{0.5\linewidth}
    \setlength{\plotheight}{0.15\paperheight}

    \pgfplotsset{
        histogram_alpha_2/.style={
                ybar interval,
                label style={font=\footnotesize},
                xmax=10,
                xtick={0,1,2,3,4,5,6,7,8,9,10},
                xticklabels={},
                x tick label style={rotate=45, anchor=east, font=\footnotesize},
                hide obscured x ticks=false,
                ymin=0,
                tick label style={font=\footnotesize},
                width=\plotwidth, height=\plotheight,
                grid=none,
            }
    }

    \begin{tikzpicture}

        \begin{axis}[
                histogram_alpha_2,
                name = left_plot,
                ylabel={Occurences, $N_h=100$},
                ymax=25,
            ]

            \addplot[color=graph_1, fill=graph_1!30!white]
            table [y={iter_100}, col sep=comma, comment chars={a}] {1D_alpha2/results_alpha2.csv};

            \draw[color=graph_4, densely dashed] (8.2, 0) -- (8.2, 25);

            \node[rectangle, draw, gray, anchor=north west, outer sep=2pt, inner sep=1pt] at (rel axis cs:0, 1)
            {
                \tikz
                \draw[color=graph_4, densely dashed] (0, 0) -- (0.5, 0)
                node[anchor=west, font=\footnotesize, color=black]{avg. gain: $1300\%$};
            };

        \end{axis}

        \begin{axis}[
                histogram_alpha_2,
                name = right_plot,
                at=(left_plot.outer east), anchor=outer west,
                ylabel={\phantom{Occurences, $N_h=100$}},
                ymax=25,
            ]

            \addplot[color=graph_2, fill=graph_2!30!white]
            table [y={cpu_100}, col sep=comma, comment chars={a}] {1D_alpha2/results_alpha2.csv};

            \draw[color=graph_4, densely dashed] (7, 0) -- (7, 25);

            \node[rectangle, draw, gray, anchor=north west, outer sep=2pt, inner sep=1pt] at (rel axis cs:0, 1)
            {
                \tikz
                \draw[color=graph_4, densely dashed] (0, 0) -- (0.5, 0)
                node[anchor=west, font=\footnotesize, color=black]{avg. gain: $500\%$};
            };

        \end{axis}

        \begin{axis}[
                histogram_alpha_2,
                name = left_plot_2,
                at=(left_plot.outer south), anchor=outer north,
                ylabel={Occurences, $N_h=200$},
                ymax=23,
            ]

            \addplot[color=graph_1, fill=graph_1!30!white]
            table [y={iter_200}, col sep=comma, comment chars={a}] {1D_alpha2/results_alpha2.csv};

            \draw[color=graph_4, densely dashed] (6, 0) -- (6, 23);

            \node[rectangle, draw, gray, anchor=north west, outer sep=2pt, inner sep=1pt] at (rel axis cs:0, 1)
            {
                \tikz
                \draw[color=graph_4, densely dashed] (0, 0) -- (0.5, 0)
                node[anchor=west, font=\footnotesize, color=black]{avg. gain: $200\%$};
            };

        \end{axis}

        \begin{axis}[
                histogram_alpha_2,
                name = right_plot_2,
                at=(right_plot.outer south), anchor=outer north,
                ylabel={\phantom{Occurences, $N_h=200$}},
                ymax=20,
                ytick={0,5,10,15},
            ]

            \addplot[color=graph_2, fill=graph_2!30!white]
            table [y={cpu_200}, col sep=comma, comment chars={a}] {1D_alpha2/results_alpha2.csv};

            \draw[color=graph_4, densely dashed] (4.789, 0) -- (4.789, 20); 

            \node[rectangle, draw, gray, anchor=north west, outer sep=2pt, inner sep=1pt] at (rel axis cs:0, 1)
            {
                \tikz
                \draw[color=graph_4, densely dashed] (0, 0) -- (0.5, 0)
                node[anchor=west, font=\footnotesize, color=black]{avg. gain: $88\%$};
            };

        \end{axis}

        \begin{axis}[
                histogram_alpha_2,
                name = left_plot_3,
                at=(left_plot_2.outer south), anchor=outer north,
                ylabel={Occurences, $N_h=400$},
                ymax=21,
            ]

            \addplot[color=graph_1, fill=graph_1!30!white]
            table [y={iter_400}, col sep=comma, comment chars={a}] {1D_alpha2/results_alpha2.csv};

            \draw[color=graph_4, densely dashed] (4.789, 0) -- (4.789, 21);

            \node[rectangle, draw, gray, anchor=north west, outer sep=2pt, inner sep=1pt] at (rel axis cs:0, 1)
            {
                \tikz
                \draw[color=graph_4, densely dashed] (0, 0) -- (0.5, 0)
                node[anchor=west, font=\footnotesize, color=black]{avg. gain: $88\%$};
            };

        \end{axis}

        \begin{axis}[
                histogram_alpha_2,
                name = right_plot_3,
                at=(right_plot_2.outer south), anchor=outer north,
                ylabel={\phantom{Occurences, $N_h=400$}},
                ymax=19,
            ]

            \addplot[color=graph_2, fill=graph_2!30!white]
            table [y={cpu_400}, col sep=comma, comment chars={a}] {1D_alpha2/results_alpha2.csv};

            \draw[color=graph_4, densely dashed] (4.677, 0) -- (4.677, 19);

            \node[rectangle, draw, gray, anchor=north west, outer sep=2pt, inner sep=1pt] at (rel axis cs:0, 1)
            {
                \tikz
                \draw[color=graph_4, densely dashed] (0, 0) -- (0.5, 0)
                node[anchor=west, font=\footnotesize, color=black]{avg. gain: $82\%$};
            };

        \end{axis}

        \begin{axis}[
                histogram_alpha_2,
                name = left_plot_4,
                at=(left_plot_3.outer south), anchor=outer north,
                xlabel={Gain in number of iterations},
                ylabel={Occurences, $N_h=600$},
                xticklabels={ 
                        {}{$-\infty\%$},
                        {}{$-50\%$},
                        {}{$0\%$},
                        {}{$25\%$},
                        {}{$50\%$},
                        {}{$100\%$},
                        {}{$200\%$},
                        {}{$500\%$},
                        {}{$1000\%$},
                        {}{$+\infty\%$},
                    },
                ymax=19,
            ]

            \addplot[color=graph_1, fill=graph_1!30!white]
            table [y={iter_600}, col sep=comma, comment chars={a}] {1D_alpha2/results_alpha2.csv};

            \draw[color=graph_4, densely dashed] (4.878, 0) -- (4.878, 19);

            \node[rectangle, draw, gray, anchor=north west, outer sep=2pt, inner sep=1pt] at (rel axis cs:0, 1)
            {
                \tikz
                \draw[color=graph_4, densely dashed] (0, 0) -- (0.5, 0)
                node[anchor=west, font=\footnotesize, color=black]{avg. gain: $93\%$};
            };

        \end{axis}

        \begin{axis}[
                histogram_alpha_2,
                name = right_plot_4,
                at=(right_plot_3.outer south), anchor=outer north,
                xlabel={Gain in CPU time},
                ylabel={\phantom{Occurences, $N_h=600$}},
                xticklabels={ 
                        {}{$-\infty\%$},
                        {}{$-50\%$},
                        {}{$0\%$},
                        {}{$25\%$},
                        {}{$50\%$},
                        {}{$100\%$},
                        {}{$200\%$},
                        {}{$500\%$},
                        {}{$1000\%$},
                        {}{$+\infty\%$},
                    },
                ymax=19,
            ]

            \addplot[color=graph_2, fill=graph_2!30!white]
            table [y={cpu_600}, col sep=comma, comment chars={a}] {1D_alpha2/results_alpha2.csv};

            \draw[color=graph_4, densely dashed] (4.861, 0) -- (4.861, 19);

            \node[rectangle, draw, gray, anchor=north west, outer sep=2pt, inner sep=1pt] at (rel axis cs:0, 1)
            {
                \tikz
                \draw[color=graph_4, densely dashed] (0, 0) -- (0.5, 0)
                node[anchor=west, font=\footnotesize, color=black]{avg. gain: $92\%$};
            };

        \end{axis}

    \end{tikzpicture}
    \caption{%
        Statistics of the gains
        $G_\text{iter}$ in number of iterations (left panels)
        and $G_\text{CPU}$ in CPU time (right panels)
        for the 1D problem \eqref{1deq} with $\alpha_0=2$.
        From top to bottom, we display the results
        for meshes with $100$, $200$, $400$ and $600$ points.
        The gains are grouped in bins of varying sizes,
        and the number of occurrences in each bin is displayed as a histogram.
    }
    \label{fig:histograms_alpha_2}
\end{figure}

Recall that training is done on meshes with $200$ and $400$ points;
validation is performed by running \cref{alg:score_function_based_on_iterations}
on meshes with varying sizes, $100$, $200$, $400$ and $600$ points.
The training time is around $15$ minutes on a cloud architecture with shared GPUs.
The tolerance for the convergence of Newton's method is fixed to $10^{-6}$.
The results are displayed on \cref{fig:histograms_alpha_2}.
Namely, in the left panels, we display the gains in number of iterations;
in the right panel, the gains in CPU time are displayed.
    {\ra%
        Note that we do not take training time into account in these CPU time gains.
        Indeed, since we are solving multiple problems with the same network,
        training time can become negligible compared to the time spent solving the PDE.
        For an in-depth study of the number of simulations needed to overcome the training time,
        we refer to \cref{tab:simulations_needed_for_profit} in \cref{sec:1D_large_alpha}.%
    }
In each case, we perform $50$ runs,
corresponding to $50$ random choices of $\Phi$ and $K$.
The gains in number of iterations are
nothing but the score $S_\text{iter}$
defined in \cref{alg:score_function_based_on_iterations},
while the gains in CPU time (denoted by $S_\text{CPU}$)
are computed in the same way
but comparing the CPU time of the two runs instead of the number of iterations.
Instead of the raw gains, we report the gain percentage:
\begin{equation}
    \label{eq:gain_as_percentages}
    G_\text{iter} = \left( S_\text{iter} - 1 \right) \times 100
    \text{\quad and \quad}
    G_\text{CPU} = \left( S_\text{CPU} - 1 \right) \times 100.
\end{equation}
{\ra%
Recall from \cref{alg:score_function_based_on_iterations}
that the number of Newton iterations is
capped at $M_\text{iter} = 2000$.
Indeed, reaching $M_\text{iter}$ iterations often implies
that Newton's method will fail to converge.
Therefore, this also puts a cap on the maximum possible
gain in iterations, which makes it harder to show the true
potential of our method, namely providing an initial guess
where Newton's method almost always converges.
}

On the top panels of \cref{fig:histograms_alpha_2}, corresponding to gains on a mesh with $100$ points, we observe that the average gain when using the initial guess predicted by the FNO is around $1300\%$ in terms of number of iterations
and $500\%$ in terms of CPU time.
In these examples, we note that the gain in the number of iterations is consistently over $100\%$.
This means that initializing Newton's method with the prediction allows it to converge in at most half the number of iterations compared to using the naive initial guess.
This gain in iterations translates into a (lower) gain in CPU time in most cases.
However, in one case out of fifty (i.e. $2.5\%$ of the time),
the total CPU time is actually larger when using the predicted initial guess.
This increase in computation time is due to the fact that
each iteration of the JFNK method is quite fast (since the mesh is coarse),
which increases the relative weight of the network call in the total CPU time.
    {\rb This may also be due to the cost of a Newton iteration being
        lower when starting from a constant guess,
        for which the linear solve is quicker.
        Our prediction is not constant,
        and thus the first Newton iterations may be more expensive
        than in the classical case where the initial guess is constant.}
We expect that, the finer the mesh,
the closer the gains in number of iterations and in CPU time will be.

Indeed, these expectations are confirmed
in the bottom three panels of \cref{fig:histograms_alpha_2},
where the results are displayed, from top to bottom,
for meshes with $200$, $400$ and $600$ points.
Our predicted initial guess enables a consistent gain
in number of iterations and in CPU time,
even on cases which were not part of the training database.
Moreover, we indeed observe that the gains in number of iterations
becomes closer to the gains in CPU time as the mesh becomes finer,
since the cost of the FNO call gradually becomes irrelevant.

An important observation is that the method is less effective
on finer meshes than on coarser ones,
with gains five to ten times lower on fine meshes.
This can be attributed to the specific way of having a good initial guess
enhances Newton's method.
Indeed, Newton's method is known to converge in two phases,
see e.g. \cite{Bre2024}.
In the second phase, Newton's method starts converging,
i.e. the residual decreases with the iterations.
This second phase only starts when the residual is small enough,
i.e., after the phase where Newton's method explores the space of solutions.
This first phase corresponds to a plateau
when graphing the residual versus the iteration number.
It should be noted that the convergence phase
is longer for harder problems, e.g. on finer meshes.
Our prediction jumpstarts the convergence phase
by providing an initial guess corresponding to a small residual.
This makes it possible to reduce, or even skip,
the plateau in the convergence curves.
This explains the smaller gains for finer meshes:
the convergence phase, which is not changed by our approach,
represents a larger part of the total number of iterations.
We expect this behavior to be reproduced for larger values of $\alpha_0$,
which also correspond to harder problems.


\subsubsection{\texorpdfstring{Results for $\alpha_0=5$ and $\alpha_0=8$}{Results for alpha0=5 and alpha0=8}}
\label{sec:1D_large_alpha}

To confirm our claims from the previous section, we now treat the same problem with larger values of $\alpha_0$.
For large values of $\alpha_0$, the nonlinear part takes precedence over the linear part in the elliptic PDE.
Therefore, the problem is harder, and the classical JFNK method will take longer to converge
(if it manages to converge at all).
In this section, we run $25$ experiments for $\alpha_0=5$ and $\alpha_0=8$, each corresponding to a random choice of $\Phi$ and $K$.
The training time is around $20$ minutes, slightly increased compared to the more linear case.

\begin{figure}[ht!]
    \centering
    \setlength{\plotwidth}{0.5\linewidth}
    \setlength{\plotheight}{0.15\paperheight}

    \pgfplotsset{
        histogram_alpha_5/.style={
                ybar interval,
                label style={font=\footnotesize},
                xmax=10,
                xtick={0,1,2,3,4,5,6,7,8,9,10},
                xticklabels={},
                x tick label style={rotate=45, anchor=east, font=\footnotesize},
                hide obscured x ticks=false,
                ymin=0,
                tick label style={font=\footnotesize},
                width=\plotwidth, height=\plotheight,
                grid=none,
            }
    }

    \begin{tikzpicture}

        \begin{axis}[
                histogram_alpha_5,
                name = left_plot,
                ylabel={Occurences, $N_h=100$},
                ymax=21,
            ]

            \addplot[color=graph_1, fill=graph_1!30!white]
            table [y={iter_100}, col sep=comma, comment chars={a}] {1D_alpha5/results_alpha5.csv};

            \draw[color=graph_4, densely dashed] (8.75, 0) -- (8.75, 21);

            \node[rectangle, draw, gray, anchor=north west, outer sep=2pt, inner sep=1pt] at (rel axis cs:0, 1)
            {
                \tikz
                \draw[color=graph_4, densely dashed] (0, 0) -- (0.5, 0)
                node[anchor=west, font=\footnotesize, color=black]{avg. gain: $2650\%$};
            };

        \end{axis}

        \begin{axis}[
                histogram_alpha_5,
                name = right_plot,
                at=(left_plot.outer east), anchor=outer west,
                ylabel={\phantom{Occurences, $N_h=100$}},
                ymax=16,
            ]

            \addplot[color=graph_2, fill=graph_2!30!white]
            table [y={cpu_100}, col sep=comma, comment chars={a}] {1D_alpha5/results_alpha5.csv};

            \draw[color=graph_4, densely dashed] (8.35, 0) -- (8.35, 16);

            \node[rectangle, draw, gray, anchor=north west, outer sep=2pt, inner sep=1pt] at (rel axis cs:0, 1)
            {
                \tikz
                \draw[color=graph_4, densely dashed] (0, 0) -- (0.5, 0)
                node[anchor=west, font=\footnotesize, color=black]{avg. gain: $1800\%$};
            };

        \end{axis}

        \begin{axis}[
                histogram_alpha_5,
                name = left_plot_2,
                at=(left_plot.outer south), anchor=outer north,
                ylabel={Occurences, $N_h=200$},
                ymax=8,
                ytick={0,2,4,6},
                yticklabels={
                        0, 2, 4, \sbox0{10} \makebox[\wd0][r]{6},
                    },
            ]

            \addplot[color=graph_1, fill=graph_1!30!white]
            table [y={iter_200}, col sep=comma, comment chars={a}] {1D_alpha5/results_alpha5.csv};

            \draw[color=graph_4, densely dashed] (6.682, 0) -- (6.682, 8);

            \node[rectangle, draw, gray, anchor=north west, outer sep=2pt, inner sep=1pt] at (rel axis cs:0, 1)
            {
                \tikz
                \draw[color=graph_4, densely dashed] (0, 0) -- (0.5, 0)
                node[anchor=west, font=\footnotesize, color=black]{avg. gain: $390\%$};
            };

        \end{axis}

        \begin{axis}[
                histogram_alpha_5,
                name = right_plot_2,
                at=(right_plot.outer south), anchor=outer north,
                ylabel={\phantom{Occurences, $N_h=200$}},
                ymax=9,
                ytick={0,2,4,6,8},
                yticklabels={
                        \sbox0{$20$}
                        \makebox[\wd0][r]{$0$},
                        \makebox[\wd0][r]{$2$},
                        \makebox[\wd0][r]{$4$},
                        \makebox[\wd0][r]{$6$},
                        \makebox[\wd0][r]{$8$},
                    },
            ]

            \addplot[color=graph_2, fill=graph_2!30!white]
            table [y={cpu_200}, col sep=comma, comment chars={a}] {1D_alpha5/results_alpha5.csv};

            \draw[color=graph_4, densely dashed] (6.125, 0) -- (6.125, 9); 

            \node[rectangle, draw, gray, anchor=north west, outer sep=2pt, inner sep=1pt] at (rel axis cs:0, 1)
            {
                \tikz
                \draw[color=graph_4, densely dashed] (0, 0) -- (0.5, 0)
                node[anchor=west, font=\footnotesize, color=black]{avg. gain: $230\%$};
            };

        \end{axis}

        \begin{axis}[
                histogram_alpha_5,
                name = left_plot_3,
                at=(left_plot_2.outer south), anchor=outer north,
                ylabel={Occurences, $N_h=400$},
                ymax=14,
                ytick={0,5,10},
            ]

            \addplot[color=graph_1, fill=graph_1!30!white]
            table [y={iter_400}, col sep=comma, comment chars={a}] {1D_alpha5/results_alpha5.csv};

            \draw[color=graph_4, densely dashed] (5.543, 0) -- (5.543, 14);

            \node[rectangle, draw, gray, anchor=north west, outer sep=2pt, inner sep=1pt] at (rel axis cs:0, 1)
            {
                \tikz
                \draw[color=graph_4, densely dashed] (0, 0) -- (0.5, 0)
                node[anchor=west, font=\footnotesize, color=black]{avg. gain: $150\%$};
            };

        \end{axis}

        \begin{axis}[
                histogram_alpha_5,
                name = right_plot_3,
                at=(right_plot_2.outer south), anchor=outer north,
                ylabel={\phantom{Occurences, $N_h=400$}},
                ymax=13,
            ]

            \addplot[color=graph_2, fill=graph_2!30!white]
            table [y={cpu_400}, col sep=comma, comment chars={a}] {1D_alpha5/results_alpha5.csv};

            \draw[color=graph_4, densely dashed] (5.543, 0) -- (5.543, 13);

            \node[rectangle, draw, gray, anchor=north west, outer sep=2pt, inner sep=1pt] at (rel axis cs:0, 1)
            {
                \tikz
                \draw[color=graph_4, densely dashed] (0, 0) -- (0.5, 0)
                node[anchor=west, font=\footnotesize, color=black]{avg. gain: $150\%$};
            };

        \end{axis}

        \begin{axis}[
                histogram_alpha_5,
                name = left_plot_4,
                at=(left_plot_3.outer south), anchor=outer north,
                xlabel={Gain in number of iterations},
                ylabel={Occurences, $N_h=600$},
                xticklabels={ 
                        {}{$-\infty\%$},
                        {}{$-50\%$},
                        {}{$0\%$},
                        {}{$25\%$},
                        {}{$50\%$},
                        {}{$100\%$},
                        {}{$200\%$},
                        {}{$500\%$},
                        {}{$1000\%$},
                        {}{$+\infty\%$},
                    },
                ymax=15,
                ytick={0,5,10},
            ]

            \addplot[color=graph_1, fill=graph_1!30!white]
            table [y={iter_600}, col sep=comma, comment chars={a}] {1D_alpha5/results_alpha5.csv};

            \draw[color=graph_4, densely dashed] (6.042, 0) -- (6.042, 15);

            \node[rectangle, draw, gray, anchor=north west, outer sep=2pt, inner sep=1pt] at (rel axis cs:0, 1)
            {
                \tikz
                \draw[color=graph_4, densely dashed] (0, 0) -- (0.5, 0)
                node[anchor=west, font=\footnotesize, color=black]{avg. gain: $210\%$};
            };

        \end{axis}

        \begin{axis}[
                histogram_alpha_5,
                name = right_plot_4,
                at=(right_plot_3.outer south), anchor=outer north,
                xlabel={Gain in CPU time},
                ylabel={\phantom{Occurences, $N_h=600$}},
                xticklabels={ 
                        {}{$-\infty\%$},
                        {}{$-50\%$},
                        {}{$0\%$},
                        {}{$25\%$},
                        {}{$50\%$},
                        {}{$100\%$},
                        {}{$200\%$},
                        {}{$500\%$},
                        {}{$1000\%$},
                        {}{$+\infty\%$},
                    },
                ymax=16,
            ]

            \addplot[color=graph_2, fill=graph_2!30!white]
            table [y={cpu_600}, col sep=comma, comment chars={a}] {1D_alpha5/results_alpha5.csv};

            \draw[color=graph_4, densely dashed] (6.084, 0) -- (6.084, 16);

            \node[rectangle, draw, gray, anchor=north west, outer sep=2pt, inner sep=1pt] at (rel axis cs:0, 1)
            {
                \tikz
                \draw[color=graph_4, densely dashed] (0, 0) -- (0.5, 0)
                node[anchor=west, font=\footnotesize, color=black]{avg. gain: $220\%$};
            };

        \end{axis}

    \end{tikzpicture}
    \caption{%
        Statistics of the gains
        $G_\text{iter}$ in number of iterations (left panels)
        and $G_\text{CPU}$ in CPU time (right panels)
        for the 1D problem \eqref{1deq} with $\alpha_0=5$.
        From top to bottom, we display the results
        for meshes with $100$, $200$, $400$ and $600$ points.
        The gains are grouped in bins of varying sizes,
        and the number of occurrences in each bin is displayed as a histogram.
    }
    \label{fig:histograms_alpha_5}
\end{figure}

The results for $\alpha_0=5$ are displayed on \cref{fig:histograms_alpha_5}.
We observe a broadly similar behavior compared to the previous case:
except in one case, the predicted initial guess
consistently outperforms the naive guess
in both iteration number and computation time.
Compared to the previous $\alpha_0=2$ case, the gains are even larger,
roughly twice as large in all cases.
As expected, since the PDE is harder to solve,
the classical JFNK method remains stuck longer in the plateau phase
for lack of a suitable initial guess.
Our prediction allows Newton's method to start with a lower residual,
which skips the plateau phase (or at least greatly reduces its duration).


These results are confirmed by the experiments with $\alpha_0=8$.
Since this case is harder to deal with,
we increase the width of the FNO to $40$ (from $30$),
and decrease the learning rate to $7.5 \times 10^{-4}$ (from $10^{-3}$).
Since the distribution of the gains is similar to the $\alpha_0=5$ case,
we do not display the histograms for the sake of brevity.
Instead, we collect the values of
the average CPU time gains in \cref{tab:gain_comparison_1D}.
We also do not report the gains in number of iterations,
since they are larger than the gains in computation time.

\begin{table}[!ht]
    \centering
    \begin{tabular}{rcccccc}
        \toprule
        mesh size    &  & $\alpha_0=2$ ($50$ examples) &  & $\alpha_0=5$ ($25$ examples) &  & $\alpha_0=8$ ($25$ examples) \\
        \cmidrule(lr){1-7}
        $100$ points &  & $+500\%$                     &  & $+1800\%$                    &  & $+5000\%$                    \\
        $200$ points &  & $+88\%$                      &  & $+230\%$                     &  & $+600\%$                     \\
        $400$ points &  & $+82\%$                      &  & $+150\%$                     &  & $+230\%$                     \\
        $600$ points &  & $+92\%$                      &  & $+220\%$                     &  & $+250\%$                     \\
        \bottomrule
    \end{tabular}
    \caption{%
        Average gains in CPU time when
        using the predicted initial guess rather than the naive one,
        for different values of $\alpha_0$.
            {\ra
                For $\alpha_0=2$ and $\alpha_0=5$,
                the hyperparameters are given in \cref{tab:hyp_final}.
                For $\alpha_0=8$,
                the width $n_p$ of the FNO is increased to $40$ (from $30$)
                and the initial learning rate $\ell_0$
                is decreased to $7.5 \times 10^{-4}$ (from $10^{-3}$).
            }
    }
    \label{tab:gain_comparison_1D}
\end{table}

In \cref{tab:gain_comparison_1D}, we observe that the mean gains in CPU time
are larger for coarser meshes and larger values of $\alpha_0$.
This is consistent with our previous observations.
Overall, in each situation under consideration,
our prediction never fails to reduce the number of iterations.
There is a small percentage (around $2\%$) of simulations where
the CPU time is increased by using our prediction.
This could be remedied by a larger network or a longer training time.



    {\ra %
        Lastly, we show that, despite the cost of training the FNO,
        the method remains cost-effective should
        multiple simulations be needed,
        especially since the inference cost is negligible.
        This is the case, for instance,
        for time-dependent problems,
        where Newton's method has to be used at each time step.
        Another use case is uncertainty quantification,
        where the PDE has to be solved for
        many realizations of the input data.
        To illustrate this, we compute the number of simulations
        needed for our method to become profitable despite the
        cost of training the FNO.
        The results are reported in
        \cref{tab:simulations_needed_for_profit}.
        We observe that, as the mesh grows finer
        and the nonlinearity increases,
        the number of simulations required for our method
        to be profitable decreases, to reach around $25$
        in the hardest cases.
        Moreover, recall that our method is resolution-invariant,
        and is therefore able to tackle multiple
        resolutions with one training.
        This is confirmed by the last row of
        \cref{tab:simulations_needed_for_profit},
        where the required number of simulations falls below $10$
        when using the FNO-based initialization on multiple grids.

        \begin{table}[!ht]
            \ra
            \centering
            \begin{tabular}{rcccccc}
                \toprule
                mesh size           &  & $\alpha_0=2$ &  & $\alpha_0=5$ &  & $\alpha_0=8$ \\
                \cmidrule(lr){1-7}
                $100$ points        &  &
                139                 &  &
                29                  &  &
                21                                                                        \\
                $200$ points        &  &
                154                 &  &
                58                  &  &
                33                                                                        \\
                $400$ points        &  &
                67                  &  &
                32                  &  &
                26                                                                        \\
                $600$ points        &  &
                44                  &  &
                27                  &  &
                26                                                                        \\
                all four mesh sizes &  &
                20                  &  &
                9                   &  &
                7                                                                         \\
                \bottomrule
            \end{tabular}
            \caption{%
                Number of simulations
                (i.e., of Newton's method calls),
                for each value of $\alpha_0$ and each mesh size,
                needed for our method to become profitable
                despite the FNO training cost.
                Since our method is able to handle multiple resolutions
                without additional training,
                the last row reports the number of simulations
                needed for our method to become profitable
                when solving the problem on all four mesh sizes
                with a single trained FNO.
            }
            \label{tab:simulations_needed_for_profit}
        \end{table}
    }

    {\rall
        \subsubsection{Comparison with other methods}
        \label{sec:comparison}

        This section is dedicated to a comparison between the naive initial guess,
        our method, and two additional methods from the literature:
        Int-Deep \cite{HuaWanYan2020} and PINN \cite{raissi2019physics}.

        Int-Deep uses a partially-trained network based on the Deep Ritz method  \cite{EYu2018}
        to provide an initial guess in Newton's method.
        The method we call ``PINN'' fully trains a PINN to approximate
        the solution to the discretized PDE
        for one specific source term and diffusion coefficient.
        This PINN is then used as an initial guess for Newton's method.
        In both cases, for each Newton solve,
        the network has to be trained on the current data
        (for \num{400} epochs for Int-Deep, and \num{10000} epochs for PINN),
        which imposes an additional, non-negligible cost.
        This cost will be higher for PINN than for Int-Deep, due to the larger number of epochs.
        However, we expect PINN to provide a better initial guess than Int-Deep,
        thus leading to fewer Newton iterations and a smaller failure rate.

        For these four methods, inference time will be negligible.
        Therefore, it is meaningful to first compare their training times and the training methods.
        This information is reported in \cref{tab:training_comparison}.
        First off, the naive method, consisting in a constant initial guess equal to $1$,
        does not require any training.
        Namely, we note that the training time for Int-Deep and PINN grows linearly
        with the number of required initial guesses to generate,
        while our method requires a single training for multiple problems.
        This behavior was already commented on when discussing
        \cref{tab:simulations_needed_for_profit}.
        Namely, it means that our method may be more suitable
        for time-dependent problems or uncertainty quantification,
        where many initial guesses are required.

        \begin{table}[!ht]
            \rall
            \centering
            \begin{tabular}{rccc}
                \toprule
                method   & training frequency         & training time            & time to generate $N$ initial guesses \\
                \cmidrule(lr){1-4}
                naive    & No training required       & \SI{0}{\s}               & \SI{0}{\s}                           \\
                Int-Deep & Once for each problem      & 400 epochs, \SI{8}{\s}   & $N \times \SI{8}{\s}$                \\
                PINN     & Once for each problem      & 10k epochs, \SI{3}{\min} & $N \times \SI{3}{\min}$              \\
                ours     & Once for multiple problems & \SI{15}{\min}            & \SI{15}{\min}                        \\
                \bottomrule
            \end{tabular}
            \caption{Comparison of training types and times for different methods.}
            \label{tab:training_comparison}
        \end{table}

        Now, we seek to compare the quality of the initial guess provided by each method.
        To do so, we run $40$ experiments with each approach,
        and count the average number of iterations $\bar N_\text{iter}$
        required for Newton's method to converge,
        as well as the failure rate
        (a simulation is said to fail if the number of iterations reaches $2000$).
        These figures are reported
        in \cref{tab:comparison_PINN_100_pts,tab:comparison_PINN_200_pts,tab:comparison_PINN_400_pts,tab:comparison_PINN_600_pts},
        for each mesh size ($100$, $200$, $400$ and $600$ points),
        and for each value of $\alpha_0$ ($2$, $5$ and $8$).
        We observe that our method consistently outperforms Int-Deep,
        probably due to the fact that Int-Deep is based on a partially-trained network.
        For a single simulation, Int-Deep will be more cost-effective than our method,
        despite the larger amount of Newton iterations required.
        Moreover, our method also outperforms PINN in most cases
        (in terms of average number of iterations),
        while its failure rate is always zero.
        As a conclusion, in a context where many initial guesses are required,
        our method is the most cost-effective.

        \begin{table}[!ht]
            \rall
            \centering
            \begin{tabular}{ccccccc}
                \toprule
                \multirow{2}{*}[-2pt]{\makecell{initial                                                                                          \\ condition}}
                         & \multicolumn{2}{c}{$\alpha_0=2$}
                         & \multicolumn{2}{c}{$\alpha_0=5$}
                         & \multicolumn{2}{c}{$\alpha_0=8$}                                                                                      \\
                \cmidrule{2-7}
                         & $\bar N_\text{iter}$             & \% failure & $\bar N_\text{iter}$ & \% failure & $\bar N_\text{iter}$ & \% failure \\
                \cmidrule(lr){1-1}\cmidrule(lr){2-3}\cmidrule(lr){4-5}\cmidrule(lr){6-7}
                naive    & 280                              & 0          & 1293                 & 30         & 1825                 & 75         \\
                Int-Deep & 276                              & 10         & 304                  & 10         & 364                  & 7.5        \\
                PINN     & 81                               & 2.5        & \bf 40               & \bf 0      & 98                   & 2.5        \\
                ours     & \bf 20                           & \bf 0      & 47                   & 0          & \bf 37               & \bf 0      \\
                \bottomrule
            \end{tabular}
            \caption{Comparison of several methods for different values of $\alpha_0$ and $N_h = 100$ points. Bold entries correspond to the best method in terms of average number of iterations or failure rate.}
            \label{tab:comparison_PINN_100_pts}
        \end{table}

        \begin{table}[!ht]
            \rall
            \centering
            \begin{tabular}{ccccccc}
                \toprule
                \multirow{2}{*}[-2pt]{\makecell{initial                                                                                          \\ condition}}
                         & \multicolumn{2}{c}{$\alpha_0=2$}
                         & \multicolumn{2}{c}{$\alpha_0=5$}
                         & \multicolumn{2}{c}{$\alpha_0=8$}                                                                                      \\
                \cmidrule{2-7}
                         & $\bar N_\text{iter}$             & \% failure & $\bar N_\text{iter}$ & \% failure & $\bar N_\text{iter}$ & \% failure \\
                \cmidrule(lr){1-1}\cmidrule(lr){2-3}\cmidrule(lr){4-5}\cmidrule(lr){6-7}
                naive    & 351                              & 0          & 781                  & 5          & 1288                 & 32.5       \\
                Int-Deep & 442                              & 12.5       & 702                  & 22.5       & 681                  & 17.5       \\
                PINN     & 198                              & 2.5        & 219                  & 0          & 306                  & 2.5        \\
                ours     & \bf 117                          & \bf 0      & \bf 159              & \bf 0      & \bf 184              & \bf 0      \\
                \bottomrule
            \end{tabular}
            \caption{Comparison of several methods for different values of $\alpha_0$ and $N_h = 200$ points. Bold entries correspond to the best method in terms of average number of iterations or failure rate.}
            \label{tab:comparison_PINN_200_pts}
        \end{table}

        \begin{table}[!ht]
            \rall
            \centering
            \begin{tabular}{ccccccc}
                \toprule
                \multirow{2}{*}[-2pt]{\makecell{initial                                                                                          \\ condition}}
                         & \multicolumn{2}{c}{$\alpha_0=2$}
                         & \multicolumn{2}{c}{$\alpha_0=5$}
                         & \multicolumn{2}{c}{$\alpha_0=8$}                                                                                      \\
                \cmidrule{2-7}
                         & $\bar N_\text{iter}$             & \% failure & $\bar N_\text{iter}$ & \% failure & $\bar N_\text{iter}$ & \% failure \\
                \cmidrule(lr){1-1}\cmidrule(lr){2-3}\cmidrule(lr){4-5}\cmidrule(lr){6-7}
                naive    & 1157                             & 5          & 1888                 & 70         & 2000                 & 100        \\
                Int-Deep & 803                              & 15         & 1117                 & 20         & 1491                 & 37.5       \\
                PINN     & \bf 495                          & 2.5        & \bf 696              & \bf 0      & 950                  & 2.5        \\
                ours     & 615                              & \bf 0      & 755                  & 0          & \bf 606              & \bf 0      \\
                \bottomrule
            \end{tabular}
            \caption{Comparison of several methods for different values of $\alpha_0$ and $N_h = 400$ points. Bold entries correspond to the best method in terms of average number of iterations or failure rate.}
            \label{tab:comparison_PINN_400_pts}
        \end{table}

        \begin{table}[!ht]
            \rall
            \centering
            \begin{tabular}{ccccccc}
                \toprule
                \multirow{2}{*}[-2pt]{\makecell{initial                                                                                          \\ condition}}
                         & \multicolumn{2}{c}{$\alpha_0=2$}
                         & \multicolumn{2}{c}{$\alpha_0=5$}
                         & \multicolumn{2}{c}{$\alpha_0=8$}                                                                                      \\
                \cmidrule{2-7}
                         & $\bar N_\text{iter}$             & \% failure & $\bar N_\text{iter}$ & \% failure & $\bar N_\text{iter}$ & \% failure \\
                \cmidrule(lr){1-1}\cmidrule(lr){2-3}\cmidrule(lr){4-5}\cmidrule(lr){6-7}
                naive    & 1736                             & 45         & 2000                 & 100        & 2000                 & 100        \\
                Int-Deep & 1159                             & 22.5       & 1614                 & 32.5       & 1866                 & 60         \\
                PINN     & \bf 769                          & 2.5        & 1125                 & 0          & 1520                 & 17.5       \\
                ours     & 899                              & \bf 0      & \bf 645              & \bf 0      & \bf 571              & \bf 0      \\
                \bottomrule
            \end{tabular}
            \caption{Comparison of several methods for different values of $\alpha_0$ and $N_h = 600$ points. Bold entries correspond to the best method in terms of average number of iterations or failure rate.}
            \label{tab:comparison_PINN_600_pts}
        \end{table}

    }

\subsection{Results in two space dimensions}
\label{sec:2D_results}

We now tackle the extension of the 1D case
\eqref{1deq} to two space dimensions,
namely the elliptic system \eqref{2deq}, with $p=2$.
It is discretized using a classical finite difference scheme
for anisotropic equations,
see the review paper \cite{EsKorBla2014}.
The test case is similar to the one we formulated in 1D.
The diffusion matrix~$K$ is obtained by
first generating a function $\delta$ as in 1D,
and then generating a random, constant symmetric positive matrix $B$;
the resulting matrix $K$ is given for all $x \in \Omega$ by $K(x)=\delta(x) B$.
This enables us to control the anisotropy of the problem
through the ratio of the eigenvalues of $B$.

Just like before, we first display the capability
of our network to predict the solution in \cref{sec:2D_prediction}.
Then, in \cref{sec:2D_initial_guess},
we compare the performance of Newton's method
with the naive initial guess and
with our predicted initial guess.
This time, the initial learning rate is set to $\ell_0 = 8 \times 10^{-4}$
to account for the increased difficulty of the problem.
Moreover, in this case, we replaced
the discretization-informed loss function \eqref{eq_loss1}
with the $H^1$ loss function \eqref{eq_loss3}.
Indeed, the value of the discretization-informed loss function
depends on the mesh size,
which prevented us from finding good hyperparameters,
even with a grid search.

\subsubsection{FNO prediction}
\label{sec:2D_prediction}

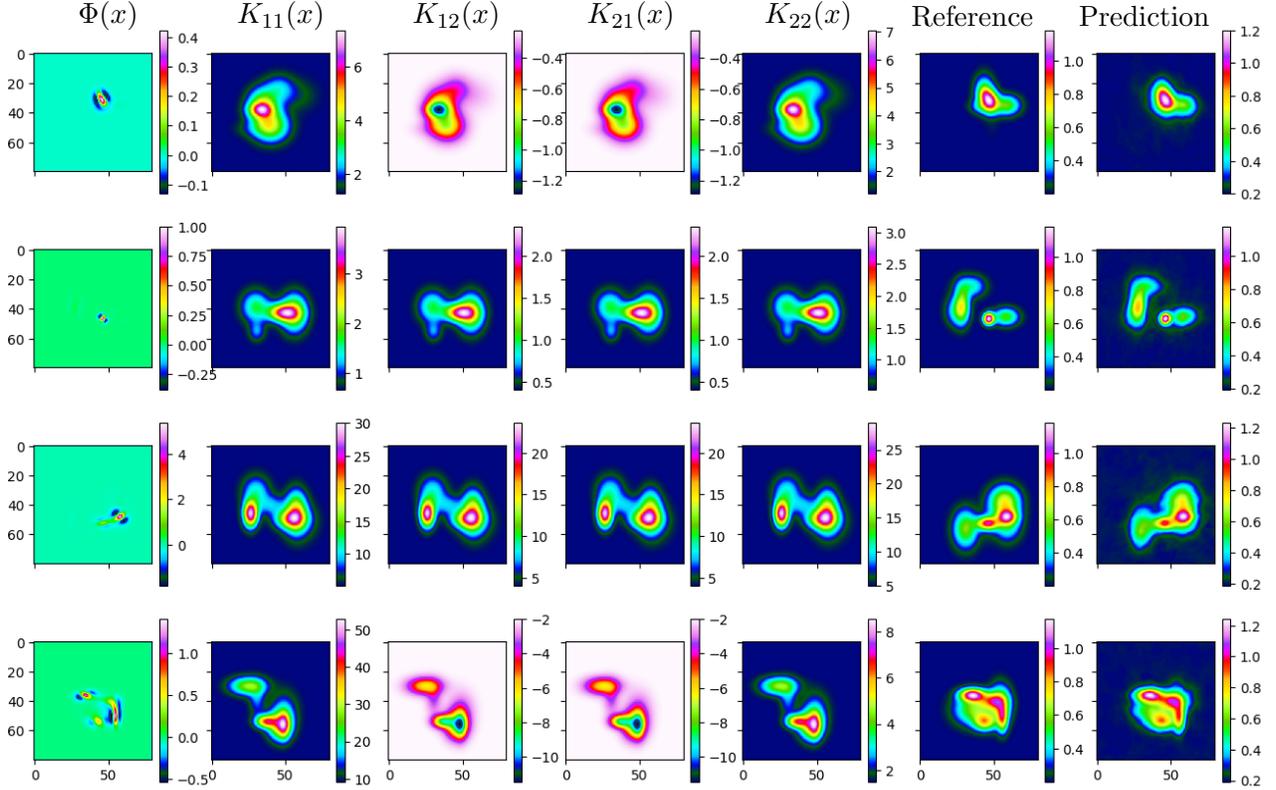
\begin{figure}[ht!]
    \centering
    \begin{tikzpicture}
        \node (fig) at (0,0) {\includegraphics[width=0.98\textwidth]{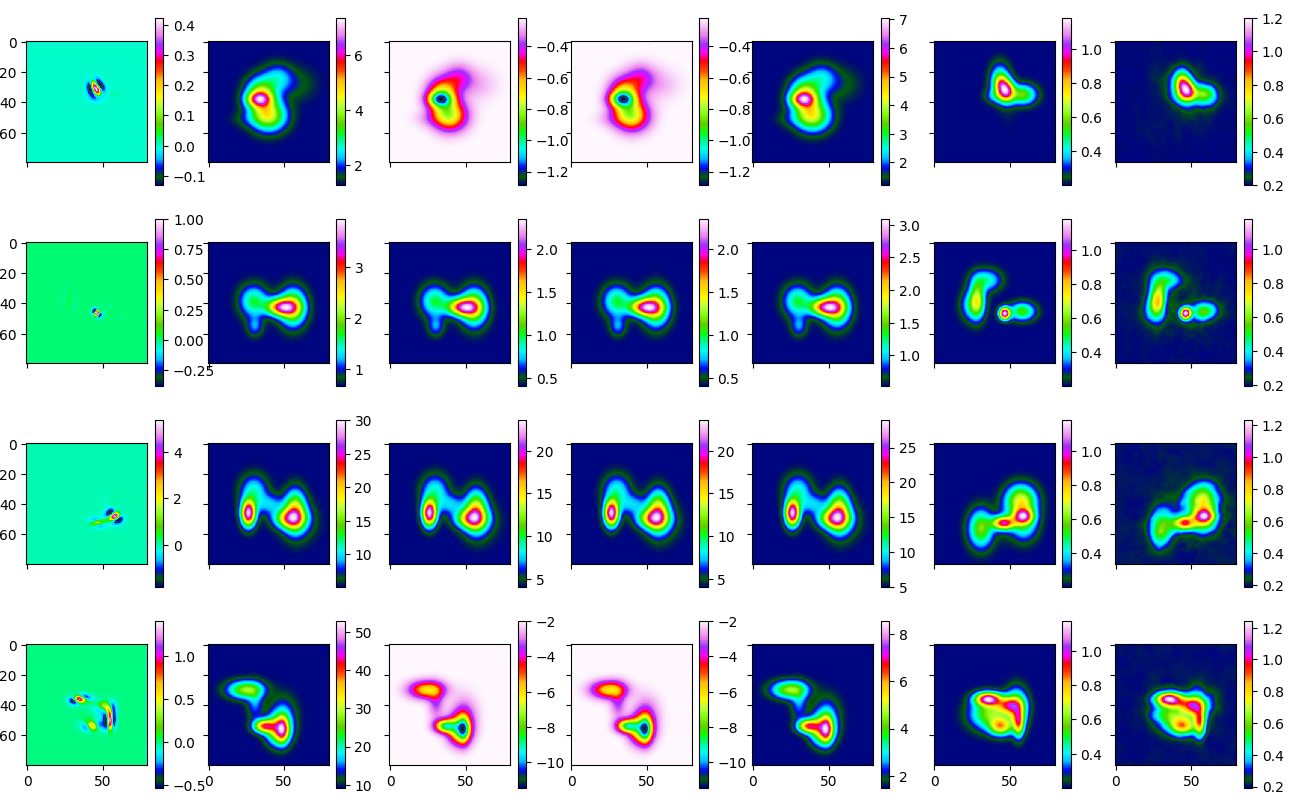}};
        \node [yshift=-5pt] at (fig.143) {$\Phi(x)$};
        \node [yshift=-5pt] at (fig.132) {$K_{11}(x)$};
        \node [yshift=-5pt] at (fig.115) {$K_{12}(x)$};
        \node [yshift=-5pt] at (fig.92) {$K_{21}(x)$};
        \node [yshift=-5pt] at (fig.68) {$K_{22}(x)$};
        \node [yshift=-5pt] at (fig.51) {Reference};
        \node [yshift=-5pt] at (fig.39) {Prediction};
    \end{tikzpicture}
    \caption{%
        From top to bottom,
        we display four random examples of generated functions for \eqref{2deq},
        as well as their predictions by the FNO.
        From left to right, we display
        the normalized values of $\Phi(x)$,
        the values of the four coordinates of the matrix $K$,
        denoted by $K_{11}(x)$, $K_{12}(x)$, $K_{21}(x)$ and $K_{22}(x)$,
        the reference solution $U(x)$,
        and the prediction of the FNO.
        In each case, we observe a good match
        between the prediction (rightmost column)
        and the reference solution (column left of the rightmost one).
    }
    \label{fig:2Dprediction}
\end{figure}

We first display, in \cref{fig:2Dprediction},
some random examples of functions $U$
predicted by the trained FNO.
They show that the FNO predicts a reasonable approximation
of the solution $U$ for variable matrices $K(x)$
(corresponding to mostly isotropic diffusion in the top three examples,
and to anisotropic diffusion in the bottom one).
We also observe that the FNO does not produce a perfect approximation
of the constant parts of the solution.

\subsubsection{Using the FNO prediction as initial guess}
\label{sec:2D_initial_guess}

Equipped with a suitable prediction,
we now study the improvement
in the convergence of Newton's method
when using it as an initial guess
rather than a naive guess.
This time, Newton's method is said to converge
if the error has reached the tolerance of $10^{-5}$.
The FNO was trained on meshes with~$60^2$ and $70^2$ points,
and the validation is performed on meshes with
$40^2$, $60^2$, $80^2$ and $100^2$ points.

In \cref{tab:2D_gains_iter_and_CPU},
we give the improvement ratios for several mesh sizes,
in terms of number of iterations and CPU times.
This table has been produced by averaging out
the results of 25 simulations.
For brevity, we do not give a detailed breakdown of the results, and we give the raw multiplicative gains $S_\text{iter}$ and $S_\text{CPU}$ rather than their percentage form \eqref{eq:gain_as_percentages}.

\begin{table}[!ht]
    \centering
    \begin{tabular}{ccccccccc} \toprule
                &  & \multicolumn{3}{c}{$S_\text{iter}$} &      & \multicolumn{3}{c}{$S_\text{CPU}$}                         \\
        \cmidrule(lr){3-5}\cmidrule(lr){7-9}
        $N_h$   &  & min                                 & avg  & max                                &  & min  & avg  & max  \\
        \cmidrule(lr){1-1}\cmidrule(lr){3-5}\cmidrule(lr){7-9}
        $40^2$  &  & 3.76                                & 4.78 & 5.86                               &  & 1.73 & 2.42 & 2.98 \\
        $60^2$  &  & 2.69                                & 3.29 & 3.95                               &  & 1.62 & 1.94 & 2.31 \\
        $80^2$  &  & 1.96                                & 2.44 & 3.02                               &  & 1.36 & 1.66 & 2.05 \\
        $100^2$ &  & 1.63                                & 2.18 & 2.85                               &  & 1.22 & 1.58 & 1.83 \\
        \bottomrule
    \end{tabular}
    \caption{%
        Total gain factor, in CPU time and number of iterations,
        obtained by using the network rather than
        the naive initial guess.
        We observe that the minimum gain
        is always greater than $1$,
        in terms of CPU time or number of iterations.
    }
    \label{tab:2D_gains_iter_and_CPU}
\end{table}

These 2D results are somewhat similar to the 1D ones,
even if an exact comparison is hard to perform
since there are more computation points in 2D than in 1D.
The main takeaway is that, in every case,
our initialization makes it possible to see
consistent gains in both numbers of iterations and CPU time.
Furthermore, it can be seen that the difference between the gain in
number of iterations and the gain in computation time
is larger than in 1D.
This is due to the two-phase convergence of
Newton's method;
the remainder of this section is dedicated to
a finer study of both these phases and
of how our prediction affects them.



\begin{table}[!ht]
    \centering
    \begin{tabular}{ccccccccccccc} \toprule
                & \multicolumn{3}{c}{$F \geq 10$} & \multicolumn{3}{c}{$1 \leq F < 10$} & \multicolumn{3}{c}{$0.1 \leq F < 1$} & \multicolumn{3}{c}{$F \leq 0.1$}                                                             \\
        \cmidrule(lr){2-4}\cmidrule(lr){5-7}\cmidrule(lr){8-10}\cmidrule(lr){11-13}
        $N_h$   & min                             & avg                                 & max                                  & min                              & avg   & max   & min  & avg   & max   & min  & avg  & max  \\
        \cmidrule(lr){1-1}\cmidrule(lr){2-4}\cmidrule(lr){5-7}\cmidrule(lr){8-10}\cmidrule(lr){11-13}
        $40^2$  & 24.50                           & 54.19                               & $\infty$                             & 15.25                            & 49.12 & 69.00 & 9.14 & 16.59 & 23.33 & 0.88 & 1.08 & 1.30 \\
        $60^2$  & 33.50                           & 65.45                               & $\infty$                             & 13.17                            & 21.72 & 33.50 & 5.93 & 9.17  & 13.60 & 0.95 & 1.11 & 1.24 \\
        $80^2$  & 16.00                           & 36.91                               & 73.00                                & 7.60                             & 12.72 & 16.60 & 3.81 & 5.83  & 8.00  & 0.89 & 1.07 & 1.28 \\
        $100^2$ & 13.57                           & 27.67                               & 69.00                                & 6.14                             & 9.85  & 15.80 & 3.31 & 4.77  & 7.08  & 0.89 & 1.12 & 1.49 \\
        \bottomrule
    \end{tabular}
    \caption{%
        Minimum, average and maximum gains
        in the number of iterations
        required to reach a certain value
        of the residual $F$ in Newton's method.
        A value of $\infty$ means that
        the predicted initial condition
        was already below the given threshold.
        We observe that the gains are mostly obtained
        for large values of the residual.
    }
    \label{tab:2D_gains_wrt_loss}
\end{table}

\begin{figure}[!ht]
    \centering
    \includegraphics{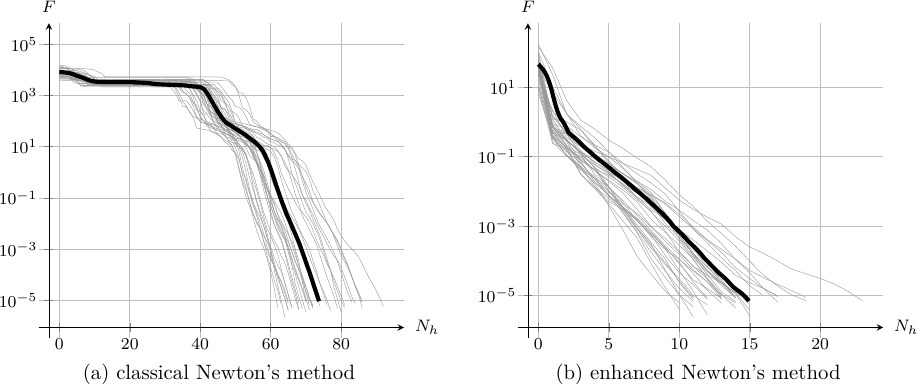}
    \caption{%
        Residual $F$ of the discretized PDE versus the number of Newton iterations
        for the classical Newton's method (left panel)
        and the enhanced Newton's method (right panel),
        for $N_h = 40^2$ points.
        The thin lines correspond to the different examples,
        while the thick line is the average over all the examples.
    }
    \label{fig:basic_vs_agent_N_40}
\end{figure}

\begin{figure}[!ht]
    \centering
    \includegraphics{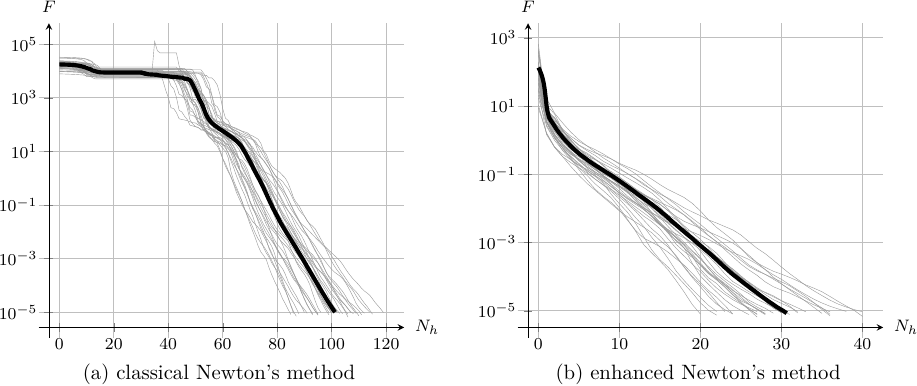}
    \caption{%
        Residual $F$ of the discretized PDE versus the number of Newton iterations
        for the classical Newton's method (left panel)
        and the enhanced Newton's method (right panel),
        for $N_h = 60^2$ points.
        The thin lines correspond to the different examples,
        while the thick line is an average over all the examples.
    }
    \label{fig:basic_vs_agent_N_60}
\end{figure}

\begin{figure}[!ht]
    \centering
    \includegraphics{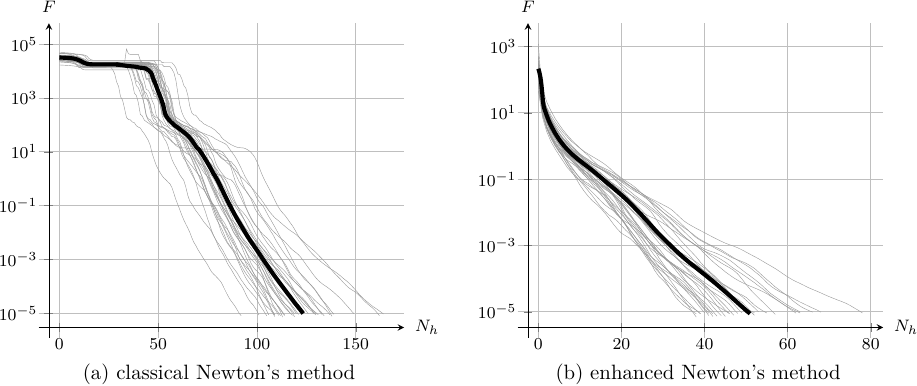}
    \caption{%
        Residual $F$ of the discretized PDE versus the number of Newton iterations
        for the classical Newton's method (left panel)
        and the enhanced Newton's method (right panel),
        for $N_h = 80^2$ points.
        The thin lines correspond to the different examples,
        while the thick line is an average over all the examples.
    }
    \label{fig:basic_vs_agent_N_80}
\end{figure}

\begin{figure}[!ht]
    \centering
    \includegraphics{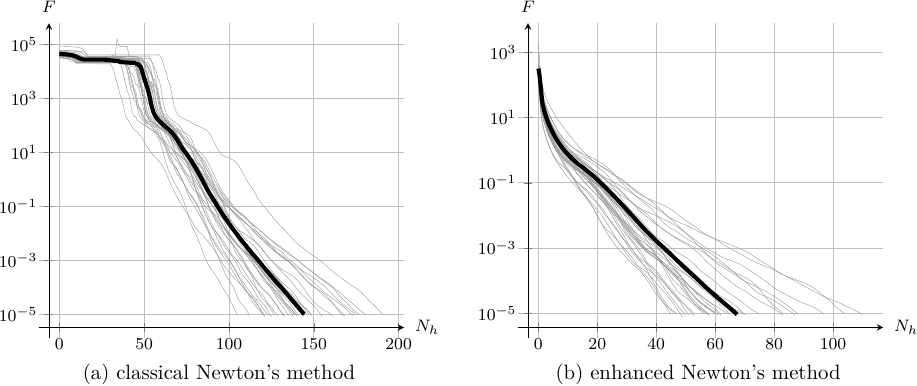}
    \caption{%
        Residual $F$ of the discretized PDE versus the number of Newton iterations
        for the classical Newton's method (left panel)
        and the enhanced Newton's method (right panel),
        for $N_h = 100^2$ points.
        The thin lines correspond to the different examples,
        while the thick line is an average over all the examples.
    }
    \label{fig:basic_vs_agent_N_100}
\end{figure}

First, in \cref{tab:2D_gains_wrt_loss},
we collect the gains in the number of iterations
required to reach some value of Newton's residual $F$.
Most of the gains are obtained in the first phase
of JFNK convergence,
i.e. to help the residual reach a small enough value
to trigger the convergence phase.
This means that, compared with a conventional initialization,
using the initial condition from the FNO successfully
avoids the first few iterations of Newton's method,
which exist to bring the residual to a small enough value.
As a consequence, once the residual is small enough,
e.g. smaller than $0.1$,
the new initialization no longer makes any difference
as Newton's method has already entered
its convergence phase.
This illustrates why the gain is lower for finer meshes.
The number of iterations below
a residual of $0.1$ is proportionally greater
as the mesh is finer.
These two phases (plateau and convergence),
and the fact that using the predicted guess bypasses
the plateau phase,
are visible on
\cref{fig:basic_vs_agent_N_40,fig:basic_vs_agent_N_60,fig:basic_vs_agent_N_80,fig:basic_vs_agent_N_100},
where we observe that reaching a residual of e.g. $10$ is much faster
when starting with the prediction rather than the naive guess.


This distribution of average gains also explains the more limited gains in computation time in~2D. Indeed, the initial iterations, when the residual is around 10, take a shorter time to complete than the subsequent ones.
This is due to the fact the solution is close to being constant, resulting in a quick linear solve (during the linear stage of the JFNK).
In contrast, the later iterations, where the residual is around 1, involve a solution further from being constant, making the associated linear solving step more challenging.
This means that, in 2D, the iterations saved by using the prediction are the fastest iterations to compute, which explains why the iteration gain is larger than the computation time gain.
The difference with 1D is that the JFNK method converges much faster in 2D once the residual reaches around $1$.
In 1D, the residual needed to be closer to $0.1$ to trigger the convergence phase.
These observations could explain why the gain in CPU time is lower than the gain in terms of iterations in 2D, even though they were comparable in 1D.


\section{Conclusion}
\label{sec:conclusion}

In this work, we tested the ability of Fourier Neural Operators (FNOs)
to predict a good initial guess for Newton's method,
applied to nonlinear elliptic partial differential equations (PDEs).
To that end, the FNO was trained to predict the PDE solution
on a large family of right-hand sides and diffusion coefficients,
in one and two space dimensions.
The solution predicted by the FNO was then used
as an initial guess in an iterative JFNK method.
To increase the accuracy of the prediction,
we used a discretization-informed loss function
(containing information about the residual of
the JFNK method applied to a discretization of the PDE)
in addition to classical data-driven loss functions.
A grid search algorithm was implemented
to select the hyperparameters.
Instead of using the prediction error
as a selection tool in the grid search,
we used a score function based on the gains,
in terms of the total number of iterations of the JFNK method,
of using the FNO prediction rather than a naive initial guess.
In all 1D and 2D test cases,
the total number of iterations decreases when
using the predicted guess,
even in the case of strong nonlinearities or anisotropy.
In the 1D test case, the worst scenario showed an average number of iterations reduced by $82\%$, as opposed to a reduction by $5000\%$ in the most favorable scenario.
The initialization mainly speeds up the plateau phase of the JFNK convergence,
rather than its convergence phase.
As a consequence, the results namely depend on the grid: the finer the grid,
the greater the number of iterations in the convergence phase, and the lesser the gains obtained with our initial guess.
In any case, the approach saves time, especially as training is quick.
The offline phase of this method (training the FNO) can be seen as the generation of a predictor, which will then be corrected online by the application of Newton's method,
to accelerate the whole process.

    {\ra%
        To further test the proposed method,
        it would be interesting to couple it with a linear or nonlinear preconditioner.
        We expect our method to still outperform the naive initialization,
        even if a preconditioner is applied.%
    }
In the future,
an interesting extension will be tackling unstructured meshes.
To that end, a possibility would be to use
Geometry-Informed Neural Operators (GINOs), recently proposed in~\cite{li2023geometry} for arbitrary geometries.
    {\rc%
        It would also be interesting to consider an application to multi-solution PDEs.
        This is a challenging problem, and our approach would require multiple databases and FNOs,
        one for each solution.
        This is in contrast to~\cite{HaoLiuYan2024},
        whose authors propose a new loss function (based on Newton's method)
        to train a neural operator, thus allowing them to tackle multi-solution PDEs.
        However, since the final approximate solution is directly given by the neural operator
        rather than by Newton's method, convergence is not ensured, contrary to our approach.%
    }
Another possible direction, on the application side, concerns time-dependent strongly nonlinear problems, such as the MHD equations in tokamaks, see \cite{franck2015energy},
or fluid flow in porous media, for which a hybrid Newton's method was designed
in \cite{LecDesFanFlaNat2023}.

\section{Conflict of interest}
\label{sec:coi}

The authors declared that they have no conflict of interest.

\bibliographystyle{plain}
\bibliography{references}

\end{document}